\documentclass[12pt]{article}
\usepackage{amsmath,amsthm,amsfonts,amssymb}
\newtheorem{thm}{Theorem}[section]
\newtheorem{lemma}[thm]{Lemma}
\newtheorem{defn}[thm]{Definition}
\newtheorem{cor}[thm]{Corollary}
\newtheorem{proposition}[thm]{Proposition}

\newtheorem{prop}[thm]{Proposition}

\theoremstyle{definition}

\newtheorem{remark}[thm]{Remark}
\theoremstyle{remark}
\DeclareMathOperator{\Fix}{Fix}
\DeclareMathOperator{\Orb}{Orb}
\DeclareMathOperator{\Diff}{Diff}
\DeclareMathOperator{\MCG}{MCG}

\DeclareMathOperator{\Homeo}{Homeo}
\DeclareMathOperator{\Int}{Int}
\DeclareMathOperator{\Per}{Per}

\newcommand{\An}{{\cal A}}

\newcommand{\R}{\mathbb R}
\newcommand{\F}{{\cal F}}

\newcommand{\cR}{{\cal R}}

\newcommand{\W}{{\cal W}}
\newcommand{\Z}{\mathbb Z}
\def\WF{\W(\F,B,C)}
\def\A{{\mathbb A}}
\def\D{{\mathbb D}}

\def\relbc{rel $B \cup C$}

\title{Fixed Points of abelian actions on $S^2$}
\author{John Franks,\thanks{Supported in part by NSF grant DMS0099640.}\ \ 
Michael Handel\thanks{Supported in part by NSF grant DMS0103435.}\ \
and Kamlesh Parwani \thanks{Supported in part by NSF grant DMS0244529.}}

\def\ti{\tilde}
\def\sinfty{S_{\infty}}
\def\gl2{GL(2, \mathbb R)}
\def\so2{SO(2)}

\def\G{{\cal G}}

\def\A{{\mathbb A}}
\def\H{{\mathbb H}}

\def\w{w}

\begin{document}
\maketitle
\begin{abstract}
We prove that if $\F$ is a finitely generated  abelian group of
orientation preserving $C^1$ diffeomorphisms of $\R^2$ which leaves
invariant a compact set then there is a common fixed point for all
elements of $\F.$ We also show that if $\F$ is any
abelian subgroup of orientation preserving $C^1$ diffeomorphisms of $S^2$ then there is a common fixed
point for all elements of a subgroup of $\F$ with  index at most two.

\end{abstract}

\section{Introduction and Notation}

The {\em global fixed point set}, denoted $\Fix(\F),$ of a group $\F$
of self maps of a space $S$ is the set of points that are fixed by
every $ f \in \F$.  Our main results are for $S = S^2$ or $S = \R^2$ and
for $\F$ a subgroup of the group $\Diff^1_+(S)$ of
orientation preserving $C^1$ diffeomorphisms of $S$.

The Lefschetz theorem implies that every element of $\Diff^1_+(S^2)$
has fixed points.  Bonatti \cite{bon:commuting} proved that if $\F$ is
a finitely generated abelian subgroup of $\Diff^1_+(S^2)$ and if a
generating set for $\F$ is sufficiently $C^1$-close to the identity,
then $\Fix(\F) \ne \emptyset$.  On the other hand, the diffeomorphisms
$R_x$ and $R_y$ of $S^2$ determined by rotation by $\pi$ about the
$x$-axis and about the $y$ -axis of $\R^3$ commute but have no common
fixed points.

If $f,g \in \Homeo_+(S^2)$ commute and $f_t, g_t$ are isotopies from
$f = f_1$ and $g = g_1$ to $id = f_0 = g_0$ define $w(f,g) \in
\pi_1(\Homeo_+(S^2),\mbox{identity}) \cong \Z/2\Z$ to be the element
determined by the loop $\gamma(t) =[f_t, g_t]$ in $\Homeo_+(S^2)$ with
basepoint at the identity.  We show below that $w(f,g)$ is independent
of the choice of the isotopies $f_t$ and $g_t$ and that if an abelian
subgroup $\F$ of $\Homeo_+(S^2)$ has a global fixed point then $w(f,g)
= 0$ for all $f,g \in \F$.

A result of the second author \cite{han:commuting} shows that if 
$f$ and $g$ are commuting orientation preserving diffeomorphisms
of $S^2$ and $w(f,g) = 0$ then $f$ and $g$ have a common fixed point.
Our main result concerning $S^2$ is a generalization to arbitrary
abelian subgroups of $\Diff^1_+(S^2)$.

\begin{thm}  \label{thm: main sphere}  For any  abelian
subgroup $\F$ of $\Diff^1_+(S^2)$ there is a subgroup $\F'$ of index
at most two such that $\Fix(\F') \ne \emptyset$. Moreover,
$\Fix(\F) \ne \emptyset$ if and only if $w(f,g) = 0$ for
all $f, g \in \F$.
\end{thm}

Druck, Fang and Firmo \cite{dff:nilpotent} generalized Bonatti's
result to include nilpotent groups.  We  conjecture that
Theorem~\ref{thm: main sphere} will also generalize to the nilpotent
case.

We derive Theorem~\ref{thm: main sphere} as a consequence of a related
result about abelian subgroups of $\Diff^1_+(R^2)$.  The Brouwer plane
translation theorem implies that if $f \in \Homeo_+(R^2)$ preserves a
compact set then $\Fix(f) \ne \emptyset$.  We generalize this as
follows.

\begin{thm}  \label{thm: partial plane} 
If  $\F$ is a finitely generated abelian
subgroup of $\Diff^1_+(\R^2)$  and  if there is a compact $\F$-invariant set $C\subset \R^2$, then  $\Fix(\F)$ is non-empty.
\end{thm}

We also prove the following corollary concerning diffeomorphisms
of closed two disk $\D^2$.

\begin{cor}  \label{cor: disk} 
If $\F$ is an abelian subgroup of $\Diff^1_+(\D^2)$
then $\Fix(\F)$ is non-empty.
\end{cor}

Related results were obtained by Hirsch \cite{hir:commuting}.

Theorem~\ref{thm: partial plane} has applications to other surfaces
because the plane is the universal cover of all surfaces but the
unpunctured sphere.  This is the reason that the proof of
Theorem~\ref{thm: main sphere} makes use of Theorem~\ref{thm: partial
plane}.  We also have the following result, which generalizes
Corollary 0.5 of \cite{han:commuting}.

\begin{thm} \label{other surfaces} Suppose that $M$ is a closed orientable surface of negative Euler characteristic, that $\F$ is a finitely generated abelian subgroup of $\Diff^1_+(M)$ and that $f \in \F$.  Let $\phi: M \to M$ be a Thurston normal form for $f$.  If there is a subsurface $M_i$ of $M$ on which $\phi$ restricts to a pseudo-Anosov homeomorphism and if $z \in \Fix(\phi|_{\Int(M_i)})$ then there exists $x \in \Fix(\F)$ such that the $f$-Nielsen class of $x$ corresponds, via the isotopy between $f$ and $\phi$, to the $\phi$-Nielsen class of $z$.
\end{thm}

Theorem~\ref{other surfaces} is in fact a special case of a result on relative Nielsen classes.  The full statement is Proposition~\ref{interior fixed point}.

We make use of the obvious fact that if $\F$ is a subgroup
of $\Homeo_+(S)$ generated by $\{f_i\}$
then $\Fix(\F) = \cap \Fix(f_i)$.  If $\F$ is generated by
$f_1,\ldots,f_n$ then we may write $\Fix(\F) = \Fix(f_1,\dots,f_n)$ and
for convenience we define $\Fix(\emptyset) = S.$

\section{Normal Form}  In this section, we assume that $S$ is a finitely punctured surface.  

We defined a normal form for an element $f \in \Diff^1_+(S)$ relative to
its fixed point set in \cite{fh:periodic}.  Given the Thurston
classification theorem \cite{Th}, the existence of a normal form for $f$
follows from the existence of a subsurface $W \subset S$ that
contains all but finitely many elements of $\Fix(f)$ and such that
$f$ is isotopic rel $\Fix(f)$ to a homeomorphism $\theta$ that
restricts to the identity on $W$.  We are now interested in abelian
subgroups of $\Diff^1_+(S)$ and it is natural to work relative to their
global fixed point sets.  Our main results are contained in
Lemma~\ref{max exists}, Lemma~\ref{max_boundary} and
Lemma~\ref{unique max} which allow us to choose $W$ in a canonical
way.

Assume  that $\F$ is a finitely generated subgroup of $\Diff^1_+(S)$ and that $B$ and $C$ are $\F$-invariant compact subsets of $S$ such that:
\begin{itemize}
\item   $B \subset \Fix(\F)$. 
\item $C \cap \Fix(\F) = \emptyset$.  
\item   $\F$ is abelian up to isotopy rel $B \cup C$. 
\end{itemize}

\begin{remark} In our applications $\F$ will be abelian but it is easier to work in the category of \lq abelian up to relative isotopy\rq.
\end{remark}

\begin{defn}       $\WF$ is the set of compact  subsurfaces   $W \subset S$ such that:
\begin{description}
\item [($\W$-1)] $\partial W$  has finitely many components, each of which is contained in and is essential in $M:=S \setminus (B \cup C)$.
\item [($\W$-2)]  $W$ contains all but finitely many points of $B$ and  every component of $W$ intersects $B$ in an infinite set.
\item [($\W$-3)] for all $f \in \F$ the following is satisfied: 
\begin{description}
\item[($\W$-3$f$)] there exists $\phi : S \to S$ isotopic to $f$ rel $B \cup C$ such that $\phi|_{W}$ is the identity.
\end{description}  
\end{description}
\end{defn}

\begin{remark} Property ($\W$-3) implies that each $W \in \WF$ is disjoint from $C$.  
\end{remark}

\begin{remark} \label{closed under intersection} If $W_i$  satisfies ($\W$-1), ($\W$-2) and ($\W$-3$f_i$) for $i=1,2$ and if    $\partial W_1$ and $\partial W_2$ intersect transversely,  then the subsurface $W$ obtained from $W_1 \cap W_2$ by removing any components that intersect $B$ in a finite set, satisfies ($\W$-1), ($\W$-2), ($\W$-3$f_1)$ and ($\W$-3$f_2$).  
\end{remark} 

Property ($\W$-2) implies that $\WF$ is empty if $B$ is finite.  The smoothness of $\F$ is used only in the following lemma.  

\begin{lemma}  If $B$ is infinite then $\WF$ is non-empty.  
\end{lemma} 

\proof  Lemma 4.1 of \cite{han:commuting} and the isotopy extension theorem (see the proof of Theorem~1.2 of \cite{fh:periodic}) imply  that  there is a subsurface $W(f)$ satisfying ($\W$-1), ($\W$-2) and ($\W$-3$f$) for any single $f \in \F$.    Since $\F$ is finitely generated, Remark~\ref{closed under intersection} produces $W \in \WF$. 
\endproof

There is a partial order on $\WF$ defined by : $W_1 < W_2$ if and only if $W_1$ is isotopic rel $B \cup C$ to a subsurface of $W_2$ but is not isotopic rel $B \cup C$   to $W_2$.   Our next lemma gives sufficient conditions for the existence of a maximal element of $\WF$.  Before stating it we recall the idea of Nielsen equivalence for fixed points of isotopic homeomorphisms.

Let $h : N \to N$ be a homeomorphism of a connected surface $N$ and
let $x,y \in \Fix(h)$.  We say that {\em $x$ is Nielsen equivalent to $y$} if there is an arc
$\alpha \subset N$ connecting $x$ to $y$ such that $h(\alpha)$ is
homotopic to $\alpha$ relative to endpoints.  Equivalently some, and
hence every, lift $\ti h : \ti N \to \ti N$ to the universal cover
that fixes a lift $\ti x$ of $x$ also fixes a lift $\ti y$ of $y$.  The Nielsen equivalence class of $x \in \Fix(f)$ is denoted $N(f,x)$.

If $\ti h$ is a lift of $h$ and $\Fix(\ti h) \ne \emptyset$ then the projection of $\Fix(\ti h)$ into $N$  is an entire Nielsen class $\mu$ of $\Fix(h)$.   We say that $\ti h$ is a {\em  lift for $\mu$} and that $\mu$ is {\em the Nielsen class determined by $\ti h$}.    Another lift  of $h$ is also a lift for $\mu$ if and only if it equals $ T \ti h T^{-1}$ for some covering translation $T$. 

Given an isotopy $h_t$  from $h_0$ to $h_1$ and a Nielsen class $\mu_0$ for $h_0$, let $\ti h$ be a lift for $\mu_0$ and let $\ti h_t$ be the lift of $h_t$ that begins with $\ti h_0$.  The  Nielsen class $\mu_t$  for $h_t$ determined by $\ti h_t$ is independent of the choice of $\ti h_0$.  We say that {\em $\mu_t$ is the Nielsen class determined by $\mu$ and $h_t$}.  If $x \in \Fix(h_s)$ is an element of $\mu_s$ and $y \in \Fix(h_t)$ is an element of $\mu_t$ then we write $N(h_s,x) \sim N(h_t,y)$.  In particular,  $N(h_t,x) \sim N(h_t,y)$ if and only if $N(h_t,x) = N(h_t,y)$.

\begin{lemma} \label{max exists} Suppose that $B$ is infinite and that either of the following is satisfied.
\begin{enumerate}
\item $C$ is finite. 
\item $\F = \langle f \rangle$,  $\Fix(f)$ is compact and $B$ is a union of Nielsen classes of $\Fix(f)$.  In particular, $B$ could be all of $\Fix(f)$.
\end{enumerate}
Then $\WF$ has maximal elements.
\end{lemma}

\proof   If $W_1 \in \WF$ is not contained in a maximal element of $\WF$ then there is an infinite increasing sequence  $W_1 \subset W_2 \subset \dots$ of non-isotopic elements of $\WF$.  We may assume that $W_l \cap B$ and the number of components of $W_l$ are independent of $l$.   The number of components of $S \setminus W_l$ is unbounded.   We may therefore choose $l$ so that some complementary component of $W_l$ is a disk $D$ that is disjoint from $B$ and from any chosen finite subset of $C$.  Since $\partial D$ is essential in $M$,  it must be that  $C$ is infinite.  

     We now assume that $\F = \langle f \rangle$ and  that   $B$ is a union of Nielsen classes of $\Fix(f)$.  Let   $C_D = C \cap D$, let $U$ be the component of $M$ that contains $\partial D$ and let $U_D= U \cup D$. In other words, $U_D$ is $U$ with the punctures in $D$ filled in.   Choose a component $\ti D$ of the full pre-image of $D$ in the universal cover $\ti U_D$ of $U_D$.  Then $\ti D$ is a disk and $C_D$ lifts to a compact subset of $\ti D$  that is invariant under a lift $\widetilde{f|_{U_D}}$   of $f|_{U_D}$.    By the Brouwer translation theorem,  there is at least one fixed point $\ti x \in \widetilde{f|_{U_D}}$.  Let $x \in U_D$ be the image of $\ti x$ and note that $x \not \in B$.  To complete the proof we will show that $N(f,x)$ contains an element of $B$ which will contradict the assumption that $B$ is a union of Nielsen classes.

There exists $\phi : S \to S$ isotopic to $f$ rel $B \cup C$ such that
$\phi|_{W_l}$ is the identity.  The isotopy from $\phi$ to $f$ lifts
to an isotopy between $\widetilde{f|_{U_D}}$ and a lift
$\widetilde{\phi|_{U_D}}$ of $\phi|_{U_D}$.  This lift setwise
preserves $\ti D$ and so has a fixed point $\ti y \in \ti D$.  Let $y
\in \Fix(\phi) \cap D$ be its projected image and let $Y$ be the
component of $W_l$ containing $\partial D$.  Choose $b \in B \cap Y$
and $z \in \partial D$.  Then
$$
N(f,b) \sim  N(\phi, b) = N(\phi,z) = N(\phi,y) \sim N(f,x)
$$ 
where the first relation follows from the fact that the isotopy
between $\phi$ and $f$ is relative to $B$, the second from the fact
that $Y \subset \Fix(\phi)$, the third from the fact that $D$ is a
$\phi$-invariant disk and the last from the definition of $\sim$.
\qed
 
\vspace{.1in}

\begin{lemma} \label{max_boundary} 
Suppose that $W$ is a maximal element of $\WF$ and that $\sigma$
is either a simple closed curve in $M$ or a simple arc whose
interior is in $M$ and whose endpoints are in $B$.
Suppose further that  the isotopy class of $\sigma$ \relbc\ is fixed by
each $f \in \F$.  Then $\sigma$ is isotopic \relbc\ to a   closed 
curve or   arc that is disjoint from $\partial W$.
\end{lemma}

\proof All isotopies in this proof are \relbc. After performing an isotopy, we may assume without
loss that $\sigma$ intersects $\partial W$ transversely and that  no component of $S \setminus (\sigma \cup \partial W)$ is a disk in the complement of $B \cup C$  whose boundary consists of  an arc in $\sigma$ and an arc in $\partial W$. We assume that  $\sigma \cap \partial W \ne  \emptyset$ and argue to a contradiction.  

Let $W'$ be the
essential subsurface obtained from a regular neighborhood of $W \cup
\sigma$ by adding in all contractible components of its complement.
Note that if $\sigma$ is an arc ending in a point of $b\in B$ not in
$W$ then $b$ is an isolated point of $B$ since all but finitely many
points of $B$ are in the interior of $W$.  Now $W \subset W'$ and $W$
is not isotopic to $W'$.  We will complete the proof by showing that
$W' \in \WF$ in contradiction to the maximality of $W$.

 Fix $f \in \F$.  Write $\sigma$ as an alternating concatenation of
subpaths $\alpha_i$ and $\beta_i$ where $\alpha_i \subset W$ and
$\beta_i \cap W = \partial \beta_i$.  By Lemma~3.5 of \cite{han:fpt} we may assume, after an
isotopy, that $f|_{\partial W \cup \sigma}$ is the identity. In
particular, $f(W) = W$.  There is a further isotopy $f_t$ from $f_0 = f$ to
$f_1 = f'$ such that $f_t(W) = W$ for all $t$ and such that $f'|_W$ is
the identity. If $C$ is a component of $W$ that contains some
$\alpha_i$ then one may assume that ${f_t}|_{\partial C}$ is the
identity for all $t$.  This follows from the fact that both $f$ and
$f'$ pointwise fix $\alpha_i$ which implies that the isotopy $f_t$ has
no net rotation about $\partial C$.  We may therefore assume that
$f_t$ pointwise fixes each $\beta_j$ for all $t$ and after an obvious
modification we may assume that $f_1|_{W'}$ is the identity.  This
proves that $W' \in \WF$ as desired.  \qed

\begin{lemma} \label{unique max} Any two maximal elements  of $\WF$ are isotopic \relbc.
\end{lemma}

\begin{proof} Suppose that $W_1$ and $W_2$ are maximal elements of $\WF$.  We
may assume by Lemma~\ref{max_boundary} that components of $\partial
W_1$ and $\partial W_2$ are disjoint or equal. Moreover, we may assume
that a component of $\partial W_1$ is isotopic \relbc\ to a 
component of $\partial W_2$ only if they are equal or bound an open 
annulus in $S \setminus (W_1 \cup W_2 \cup B \cup C).$

If $W_1 \ne W_2$ then the interior of one of them contains points not
in the interior of the other. We will show this leads to a
contradiction.  So assume without loss of generality that there is a
non-empty component of $int(W_2) \setminus W_1$ and let $Y$
denote its closure.   Then $Y$ is a compact surface with boundary, is disjoint from $C$ and contains at most finitely many elements of $B$.

Let $W_2^0$ denote the component of $W_2$ which contains $Y$.  Since
$W_2^0 \cap B$ is infinite, at least one component of $W_2^0 \setminus
Y$ must intersect $B$ in an infinite set.  Call the closure of such a 
component $Z$ and note that $Y \cap Z \subset \partial W_1$.   If $Y \cap B \ne \emptyset$ then there is an arc
$\sigma \subset W_2$ connecting an element of $B$ in $Y$ to an element
of $B$ in $Z \cap W_1$.  This contradicts Lemma~\ref{max_boundary} and
the fact that $\sigma$ is isotopic \relbc\ to its image under any
element of $\F.$ We conclude that $B \cap Y= \emptyset$.
  A similar
argument shows that $Z$ is the only component of $W_2^0 \setminus
Y$ which intersects $B$, since otherwise there would be an arc
joining points of $B$ in different components of $W_2^0 \setminus
Y$ and crossing a boundary component of $Y$ transversely in a single
point.  This would also contradict Lemma~\ref{max_boundary}.

We next assume $Y$ is not an annulus.   There
must be at least one common boundary component of $Z$ and $Y$. If
there is more than one then there is a simple closed curve $\sigma$ in
$int(Y \cup Z) \subset int(W_2^0)$ which has non trivial
intersection number with two of these common components.  Since
$\sigma \subset W_2^0$ it is isotopic \relbc\ to its image under any
element of $\F.$ But it cannot be isotoped to be disjoint from $Y \cap Z$
which contradicts Lemma~\ref{max_boundary}.

We may therefore assume that 
$Z \cap Y$ consists of a single common boundary component; call it $X$.
There is then a simple arc $\sigma$ in $int(Y \cup Z) \subset int(W_2^0)$ 
joining points of $B$, crossing $X$ transversely twice and
representing a non-trivial element of $H_1(Z \cup Y, Z)$.
We now have compact surfaces $Y \subset W_2^0 \subset S$ and 
we wish to lift to the universal covering space $\ti S.$  More
precisely, choose a component $\ti Y$ of the complete lift of $Y$
and a component $\ti W$ of the complete lift of $W_2^0$ such that
$\ti Y \subset \ti W \subset \ti S.$ Note that $\ti Y$ separates
$\ti W$ and the fact that $\sigma$ represents a non-trivial element
of $H_1(Z \cup Y,  Z)$ means that it has a lift $\ti \sigma$ to a 
simple arc in $\ti W$ with endpoints separated by $\ti Y.$  
If $\sigma'$ is a simple arc in $W_2^0$ isotopic in $S$ 
to $\sigma$\ \relbc \ then
it has a lift $\ti \sigma'$ with the same endpoints as $\ti \sigma.$
It follows that $\ti \sigma' \cap \ti Y \ne \emptyset$ and
$\sigma' \cap Y \ne \emptyset$ contradicting Lemma~\ref{max_boundary}.
This completes the proof when $Y$ is not an annulus.

If $Y$ is an annulus then each of its boundary components is a
boundary component of $W_1, \ W_2$ or both.  If one is in $\partial
W_1$ and the other is in $\partial W_2$ we would contradict the
assertion above that such boundary components can be isotopic \relbc\
only if the open annulus they bound lies in $S \setminus (W_1 \cup W_2
\cup B \cup C).$ If both components of $\partial Y$ are contained in
$\partial W_2$ then $Y$ would be a component of $W_2$ which is a
contradiction.  Hence both components of $\partial Y$ must be in
$\partial W_1$ and disjoint from $\partial W_2.$ In this case $W_2^0
\setminus int(Y)$ can have no component which is an unpunctured annulus.

Let $Z'$ be a component of  $W_2^0 \setminus int(Y)$ which is different
from $Z$.  It is disjoint from $B$ and is not an annulus.  In this case too there is an
arc $\sigma$  with endpoints in $B$ whose interior lies in 
$W_2^0 \setminus B$ which crosses $X$ transversely twice and
represents a non-trivial element of $H_1(Z \cup Y \cup Z', Z)$.
The same covering space argument given above 
shows if $\sigma'$ is a simple arc in $W_2^0$ isotopic in $S$ 
to $\sigma$\ \relbc \ then $\sigma' \cap X \ne \emptyset.$
This again contradicts Lemma~\ref{max_boundary}.
\end{proof}

The following result is proved in \cite{FLP} (see Th\'eor\`eme III of 
Expos\'e 12) at least in the case that $S$ has no punctures.  The 
proof we give here is quite different and very much in the spirit
of extracting information from the Nielsen classes of fixed points
of iterates which we continue in later sections.  More details about
Nielsen classes (including definitions) 
are given at the beginning of Section~\ref{sec:nielsen} and in references
cited there.  For the following result we use only the fact that
interior fixed points of pseudo-Anosov homeomorphisms are unique
in their Nielsen class.

\begin{lemma}  \label{pA-conj} 
Let $S$ be a surface, perhaps finitely punctured, with negative Euler
characteristic and let $f$ be a pseudo-Anosov homeomorphism of $S$. If
$g$ is another pseudo-Anosov homeomorphism of $S$ which is homotopic
to $f$ then there exists a unique homeomorphism $h: S \to S$ which is
homotopic to the identity and satisfies $h \circ f = g \circ h.$
\end{lemma} 
\begin{proof}
We consider the Poincar\'e disk model for the hyperbolic plane $\H$.  In
this model, $\H$ is identified with the interior of the unit disk and
geodesics are segments of Euclidean circles and straight lines that
meet the boundary in right angles. A choice of complete hyperbolic structure on
$S$ provides an identification of the universal cover $\ti S$ of $S$
with $\H$.  Under this identification covering translations become
isometries of $\H$ and geodesics in $S$ lift to geodesics in $\H$.  The
compactification of the interior of the unit disk by the unit circle
induces a compactification of $\H$ by the \lq circle at infinity\rq\
$\sinfty$.  Geodesics in $\H$ have unique endpoints on $\sinfty$.
Conversely, any pair of distinct points on $\sinfty$ are the endpoints
of a unique geodesic.  

Suppose that $F: S \to S$ is an orientation preserving homeomorphism
of $S.$ We can use the identification of $\H$ with $\ti M$ and write
$\ti F : \H \to \H$ for lifts of $F: S \to S$ to the universal cover.
A fundamental result of Nielsen theory is that every lift $\ti F : \H
\to \H$ extends uniquely to a homeomorphism (also called) $\ti F : \H
\cup \sinfty \to \H \cup \sinfty$ and the restriction of this
homeomorphism to $\sinfty$ depends only on the homotopy class of $F$
and the choice of lift $\ti F$.  Using the identification above we
will typically consider $\ti F$ to be a homeomorphism of $\ti S \cup
\sinfty.$ The homeomorphism $F$ induces an isomorphism (which we will
denote by $F_\#$) of the fundamental group of $S$ or equivalently of
the group of covering translations of $\ti S.$

Let $X$ denote the compactified universal covering space of $S$, i.e
$X = \ti S \cup \sinfty.$ Let $\ti f: X \to X$ be the extension of a
lift of $f$ and let $\ti g: X \to X$ be the extension of the lift of
$g$ obtained by lifting the homotopy from $f$ to $g$ starting at $\ti
f.$ If $x$ is a point of period $p$ of $f$ and $\ti x$ is a lift of $x$ then
there is a unique covering translation $T_0$ such that 
$T_0 \circ \ti f^p( \ti x) = \ti x.$ There is a unique 
point $\ti y \in \ti S$ which
satisfies $T_0 \circ \ti g^p( \ti y) = \ti y.$ This is because $f^p$ and
$g^p$ are homotopic pseudo-Anosov homeomorphisms
and for any pseudo-Anosov homeomorphism any interior fixed point is
the unique point in an essential Nielsen class.

The assignment $\ti x \mapsto \ti y$ defines a bijection
$\ti h : \widetilde \Per(f) \to \widetilde \Per(g)$ from 
the lifts of all periodic points of $f$ to the lifts of
all periodic points of $g$.  This function satisfies
$\ti h \circ \ti f = \ti g \circ \ti h$ because
\[
f^{-1}_\#(T_0) \circ \ti f^p( \ti f(\ti x)) = 
f^{-1}_\#(T_0)\circ \ti f \circ \ti f^p( \ti x) = 
\ti f \circ T_0 \circ \ti f^p( \ti x) = \ti f(\ti x).
\]
so $\ti f(\ti x)$ is the unique fixed point in $\ti S$ 
of $f^{-1}_\#(T_0) \circ \ti f^p.$
Similarly $\ti g(\ti y)$ is the unique fixed point of $g^{-1}_\#(T_0)
\circ \ti g^p$.  Since $g_\# = f_\#$ we conclude that $\ti h( \ti f(\ti x))
= \ti g(\ti y) = \ti g( \ti h(\ti x)).$  A similar computation shows
that $T \circ \ti h = \ti h \circ T$ for any covering
translation $T$.

Let $\F^u(f)$ and $\F^s(f)$ denote the $f$-invariant expanding and
contracting foliations and let $\ti \F^u(f)$ and $\ti \F^s(f)$ denote their
lifts to the universal covering space $\ti S$.  The leaves of these
foliations have well defined ends in $\sinfty.$ (See \cite{FLP} for
more details on the properties of these foliations.)
Note that any leaf of $\ti \F^u(f)$ intersects a
leaf of $\ti \F^s(f)$ in at most one point and there is a point of
intersection precisely if the ends of the two leaves are linked in
$\sinfty.$ By linked we mean that the ends (there may be more than
two) of one leaf separate the ends of the other.  Note also that
the ends of the leaf $\ti W^u(\ti x,  f)$ 
containing $\ti x$ (the fixed point of $T_0 \circ
\ti f^p$) are precisely the attracting fixed points of the
restriction of $T_0 \circ \ti f^p$ to $\sinfty$  and
hence they are also the ends of the leaf
$\ti W^u(\ti h(\ti x), g).$  This implies that linking properties
of the ends of leaves through points of $\widetilde \Per(f)$ are
preserved by $\ti h$ as are the order properties of those leaves.
By this we mean if such a leaf in $\F^u(f)$ separates two others
then the corresponding leaves in $\F^u(g)$ will have the corresponding
separation properties.

{From} this it is straightforward to see that $\ti h$ may
be extended uniquely to a bijection from a dense set of lifted heteroclinic
points in $\ti W^u(\ti x, f)$ to  lifted heteroclinic points of 
$\ti W^u(\ti h(\ti x), g)$ by defining 
$\ti h( \ti W^u(\ti x, f) \cap \ti W^s(\ti z, f)) = \ti W^u(\ti h(\ti x), g)
\cap \ti W^s(\ti h(\ti z) g)$ for all $\ti z \in \widetilde \Per(f)$ for
which $\ti W^u(\ti x, f) \cap \ti W^s(\ti z, f)$ is non-empty.

Moreover, the function $\ti h$ is order
preserving on the dense subset of $\ti W^u(\ti x, f)$ on which it
is defined and hence extends uniquely to a continuous order 
preserving function $\ti h: \ti W^u(\ti x, f) \to \ti W^u(\ti h(\ti x), g).$
Doing the same construction for all $\ti z \in \widetilde \Per(f)$
we extend $\ti h$ to the union of all unstable leaves through
points of $\widetilde \Per(f)$.  
It is clear that $\ti h$ so defined 
respects those stable foliation leaves which contain points
of $\widetilde \Per(f)$ and the restriction of $\ti h$ to the
lifted heteroclinic points in a single stable leaf is order 
preserving.  Hence there is a unique extension of $\ti h$ to the
union of all stable and unstable leaves which contain a point
of $\widetilde \Per(f)$.  Finally the local product structure
of the foliations is preserved and hence there is a unique
continuous extension of $\ti h$ to $\ti h : \ti S \to \ti S.$

If $x \in \widetilde \Per(f)$ the ends of the leaf $\ti W^u(\ti x, f)$
coincide with the ends of the leaf $\ti W^u(\ti h(\ti x), g)$ since
both are the attracting fixed points of $\ti f = \ti g$ on $\sinfty$.
Thus setting $\ti h = id$ on $\sinfty$ extends it continuously to all
of $X.$ We have $\ti h \circ \ti f = \ti g \circ \ti h$ and $T \circ
\ti h = \ti h \circ T$ for every covering translation $T$, because
these equations hold on the dense subset $\widetilde \Per(f)$.

Thus $\ti h$ is the lift of a homeomorphism $h : S \to S.$
which satisfies  $h \circ f = g \circ h$ and which is
isotopic to the identity since $\ti h$ commutes with
covering translations.  If $k$ were another such conjugacy from $f$ to $g$
then its identity lift $\ti k$ would have to agree with $\ti h$
on the dense set $\widetilde \Per(f)$ and hence they would agree 
everywhere.
\end{proof}

The following result in the case $n=2$ is essentially contained
in Lemmas 2.2 and 2.3 of \cite{han:commuting}.

\begin{lemma}  \label{virt_cyclic} 
Let $S$ be a surface, perhaps finitely punctured, with negative Euler
characteristic and $A = \langle \beta_1, \dots, \beta_n \rangle$ an abelian
subgroup of $\MCG^+(S).$ Suppose that either some element of $A$ is
pseudo-Anosov, or $A$ is finite.  Then there is a
homomorphism $\phi: A \to \Homeo^+(S)$ such that $\phi(g)$ is a
canonical Thurston representative in the mapping class $g.$
Moreover if $A$ contains a pseudo-Anosov element then
$\phi(A)$ is the direct sum of an infinite cyclic group with
a finite group.
\end{lemma}

\begin{proof}
In the case that $A = \langle \beta_1, \dots, \beta_n \rangle$ is finite,  
Kerchoff's solution of the Nielsen realization problem(\cite{Ker})
 implies that there is a Riemannian metric on $S$ and homomorphism
$\phi$ from $A$ to the isometries of this metric such that 
for each $\beta \in A,\ \phi(\beta)$ is in the mapping class $\beta$.

If $A$ contains an irreducible pseudo-Anosov element $\beta_0$, we will prove the
result when $A$ is the centralizer of $\beta_0$. Let $f_0$ be a pseudo-Anosov diffeomorphism with $\beta_0$ equal to the isotopy class $[f_0]$ of $f_0$. Suppose $\beta$ is in the centralizer of $\beta_0$ and $h$ is
a diffeomorphism with $\beta = [h].$ Then $\beta_0 = [h^{-1} \circ f_0
\circ h]$ so Lemma~\ref{pA-conj} implies there is a unique
homeomorphism $k: S \to S$, homotopic to the identity, such that
\[
(h\circ k)^{-1} \circ f_0 \circ (h \circ k) = f_0.
\]
In other words there is a unique choice of $h$ representing $\beta$ that commutes with $f_0$.   We define $\phi: A \to \Homeo^+(S)$ by letting $\phi(\beta)$ be this unique $h$.  The uniqueness property implies that $\phi$
is a homomorphism.  It is clear that $\beta = [\phi(\beta)].$

To see that $h = \phi(\beta)$ is a Thurston canonical form we observe
that it preserves both the stable and unstable foliations of $f_0$.
These foliations possess transverse invariant measures $\mu_s$ and
$\mu_u$ which are unique up to a scalar multiple (See Expos\'e 12
of  Fathi, Laudenbach,
and Poenaru \cite{FLP}).  Since $h$ carries
one such measure to another there is a unique positive constant
$\alpha_u(h)$ such that $h_*(\mu_u) = \alpha_u(h) \mu_u.$
Similarly there is a unique positive constant
$\alpha_s(h)$ such that $h_*(\mu_s) = \alpha_s(h) \mu_s.$  
The fact that $h$ is homeomorphism implies that the area of
$h(S)$ equals the area of $S$ so  $\alpha_s(h) \alpha_u(h) = 1.$

The functions $\alpha_s$ and $\alpha_u$ are homomorphisms to $\R^+,$ the
positive reals under multiplication.  If $\alpha_u(h) =1$ and $h$ has
an interior fixed point and fixes a branch of its unstable manifold
then it is clear that $h = id$ since this branch is dense. Since any
$h$ must permute fixed points of $f_0$ and the branches of their
unstable manifolds we must have $h^k = id$ for some $k >0.$ Thus the
kernel of $\alpha_u$ is contained in the torsion subgroup of $\phi(A)$
and in fact equal to it since the image of $\alpha_u$ has no torsion.
If $\alpha_u(h) \ne 1$ then $h$ is pseudo-Anosov with the same
foliations as $f_0$ (perhaps with stable and unstable switched).  To
complete the proof we need only the fact that $\alpha_u(\phi(A))$ is
cyclic.  But this follows from the fact that the image of $\alpha_u$
must be discrete.  If this were not the case there would be
pseudo-Anosov homeomorphisms with expanding and contracting constants
arbitrarily close to one.  Such a homeomorphism which preserved
branches of stable and unstable leaves at an interior fixed point would be
$C^0$ close to the identity.  This is not possible since no
pseudo-Anosov homeomorphism can be homotopic to the identity.
\end{proof}

We now apply these results to produce normal forms for an abelian
subgroup $\F$ of $\Diff^1_+(S)$ assuming that $\Fix(\F) = \emptyset$.
Our applications in this paper are in genus 0 and we use normal forms
in proofs by contradiction, ultimately showing that $\Fix(\F) \ne
\emptyset$.  Thus the proposition below is less important for future
applications than the preceding results on $\WF$, which allow us to
produce normal forms.

\begin{prop} \label{normal form}  Suppose that $S$ is a finitely punctured surface,   that $ f_1,\dots,f_n  \in \Diff^1_+(S)$ generate an abelian subgroup $\F$ and that $\Fix(\F) = \emptyset$.  Suppose further that:
\begin{itemize}
\item  $K \subset \Fix(f_1,\dots,f_{n-1})$ is compact and $\F$-invariant.
\item  $L \subset \Fix(f_n)$ is compact and  $\F$-invariant.
\item  If $L$ is infinite then $\Fix(f_n)$  is compact and $L = \Fix(f_n)$.
\item  If $(K \cup L)$ is finite then $\chi(M) < 0$ where $M:= S\setminus (K \cup L)$.     
\end{itemize}
Then there is a finite set $R$ of disjoint simple closed curves  (called {\em reducing curves}) in $M$     and  for $1 \le j \le n$ there are homeomorphisms $\theta_j : S \to S$ isotopic to $f_j$ rel $K \cup L$   
such that:
\begin{description}
\item [(1)] $\theta_j$ permutes disjoint open annulus neighborhoods  in $S$ of the elements of $R $ for   $ 1\le j \le n$.

\end{description}
Denote the union of the annular neighborhoods by $\An$, let $\{S_i\}$ be the components of $S \setminus \An$, let $X_i =
(K \cup L) \cap S_i$ and let $M_i := S_i \setminus X_i$.
\begin{description}
\item [(2)] If $X_i$ is infinite then  either  $X_i \subset L$  and $\theta_n|_{S_i}$ is the identity or $X_i \subset K$ and $\theta_j|_{S_i}$ is the identity for $ 1\le j \le n-1$.
\item [(3)] If $X_i$ is finite then $M_i$ has negative Euler characteristic  and, for $ 1\le j \le n$,  $\theta_j^{r_{ij}}|_{M_i}$  is either periodic or pseudo-Anosov, where    $r_{ij}$ is the smallest positive integer such that $\theta_j^{r_{ij}}(M_i) =  M_i$.  Moreover, the $\theta_j^{r_{ij}}|_{M_i}$'s generate an abelian subgroup that is either finite or virtually cyclic.
\end{description}

\end{prop}

\proof  During this proof $1 \le j \le n-1$.  

If $L$ is finite, define $W_L= \emptyset$.  Otherwise,  let $W_L$ be a maximal element of $\W(\langle f_n\rangle,L,K)$.  By definition, there is a diffeomorphism $\theta_n :S \to S$ that is isotopic to $f_n$ rel $K \cup L$ and that restricts to the identity on $W_L$.  Since $f_j$ preserves $L$ and $K$ and commutes with $f_n$ it follows that  $f_j(W_L)$ is a maximal element of $\W(\langle f_n\rangle,L,K)$.  Corollary~\ref{unique max} therefore implies that $f_j$ preserves $W_L$ up to isotopy rel $K \cup L$.  Choose diffeomorphisms $\theta_j:S \to S$ that are isotopic to $f_j$ rel $K \cup L$ and that preserve $W_L$.   

Denote the finitely punctured subsurface $S \setminus W_L$ by $S'$ and the subset of $L$ that is not contained in $W_L$ by $L'$.    If $K$ is finite, define $W_K = \emptyset$.  Otherwise, define $W_K$ to be a maximal element of $\W(\langle {\theta_1}_{|S'},\dots,{\theta_{n-1}}_{|S'}\rangle, K,L')$.   After another isotopy rel $K \cup L$ we may assume that each $\theta_j$ restricts to the identity on $W_K$ and that $W_K$ is $\theta_n$-invariant. 

A partial  set of reducing curves $R'$ is defined from $\partial W_K \cup \partial W_L$ by removing any pair of curves that cobound an unpunctured annulus and replacing them with a core curve of that annulus.  We may assume without loss that (1) and (2) are satisfied with respect to $R'$.   Let $S''$ be the complement in $S$ of  $W_L \cup W_K$ and the annuli associated to $R'$.  Then $L
\cup K$ intersects $S''$ in a finite set and its complement in $S''$
is denoted $N$.  To complete the proof we must show that $R'$ can be
extended by adding in simple closed curves in $N$ to maintain (1) and to arrange (3).  The existence of reducing curves for the induced action of $\F$ on $N$ is well known; see for example Lemma~2.2 of \cite{han:commuting}.  These new curves divide $N$ into irreducible subsurfaces and (3) is then a consequence of Lemma~\ref{virt_cyclic}.
\endproof

The set $R$ of Proposition~\ref{normal form} is not uniquely determined.  For example, if $\theta_j|_{M_i}$ is the identity for all $1 \le j \le n-1$ and if $\theta_n$ preserves an essential non-peripheral  simple closed curve $\alpha \subset M_i$ then one can add $\gamma$ to $R$.   Lemma~\ref{augmented regular nbhd} below  produces such curves $\alpha$ and is applied in the proof of Lemma~\ref{lem: AB=>B}.

Suppose   that  $M$ is obtained from a closed surface by puncturing at a possibly infinite set. Choose a complete hyperbolic metric of finite volume for $M$.  It is not strictly necessary to work with such a metric but it simplifies statements and proofs because the free conjugacy class of any essential non-peripheral closed curve is represented by a unique closed geodesic. In particular, a homeomorphism $f : M \to M$ induces a bijection $f_\#$ of the set of closed geodesics in $M$.

 For any finite collection $\Sigma = \{\sigma_1,\dots,\sigma_l\}$ of closed geodesics in $M$  define the {\em augmented  regular neighborhood} $N(\Sigma)$   to be the union of a regular neighborhood of $\sigma_1 \cup \dots \cup \sigma_l$  with   all contractible components of its complement.    The isotopy class of $N( \Sigma)$ is well defined.

The following lemma is  similar to Lemma~2.2 of \cite{ht:nielsen}.  In this paper we apply it only with $l = 1$ and with $M$ having genus zero.

\begin{lemma} \label{augmented regular nbhd} If $f : M \to M$ is an orientation preserving homeomorphism that preserves the homotopy class of   $\sigma_i$ for $1 \le i \le l$ then  $f$ is isotopic to  a homeomorphism $F : M \to M$ that preserves $N( \Sigma)$   and such that $F|_{N(\Sigma)}$ has finite order.
\end{lemma}

\proof  There is no loss in assuming that $\sigma_1 \cup \dots \cup \sigma_l$, and hence  $N(\Sigma)$, is connected.     The lemma is clear if $N(\Sigma)$ is an annulus so we may assume that  each element of $\Sigma$ is non-peripheral in $N(\Sigma)$.  

Define $\Sigma^c$ to be the set of simple closed geodesics $\gamma \subset M$ that are disjoint from each $\sigma_i$.  Then $\Sigma^c$ is $f_\#$-invariant.  We claim that, up to isotopy,  $N(\Sigma)$ is the unique  essential subsurface $N \subset M$ that has  finite type,  that is a closed subset of $M$ and that satisfies 
\begin{description}
\item[($\ast$)] $\gamma \in \Sigma^c$ if and only if  there is a representative of the isotopy class of $N$ that is disjoint from $\gamma$.
\end{description}

To see that $N(\Sigma)$ satisfies $(\ast)$, suppose that $\gamma \in \Sigma^c$.   Then $\gamma$  is contained in a non-contractible component of the complement of $\sigma_1 \cup \dots \cup \sigma_l$.  For an appropriate choice of regular neighborhood of $\sigma_1 \cup \dots \cup \sigma_l$, the augmented regular neighborhood is disjoint from $\gamma$.   Conversely, if  there is a representative of the isotopy class of $N$ that is disjoint from $\gamma$
 then $\gamma$ is disjoint from a representative of the isotopy class of each $\sigma_i$.  Since $\gamma$ and $\sigma_i$ are geodesics, they must be disjoint or equal.  Since $N(\Sigma)$ is not an annulus, $\gamma$ must be disjoint from each $\sigma_i$. 

To prove uniqueness up to isotopy, suppose that $N_1$ and $N_2$ satisfy ($\ast$).  After an isotopy, we may assume without loss that  $\partial N_1 \cap N_2 = \emptyset$.  Since $N_1 \cap N_2 \ne \emptyset$ (because $N_2$ must intersect any essential non-peripheral simple closed curve in $N_1$)  it follows that   $N_1 \subset N_2$.  The symmetric argument shows that $N_2$ is isotopic into $N_1$ and so $N_1 = N_2$.   

This completes the proof of the claim.  Since $N(\Sigma)$ and $f(N(\Sigma))$ both satisfy ($\ast$), $f$ is isotopic to a homeomorphism $F$ that preserves $N(\Sigma)$.

It remains to prove that the isotopy class of $F|_{N(\Sigma)}$ has
finite order.  Since $\sigma_i$ is $f_\#$-invariant, $\theta$ can not
be pseudo-Anosov and $\sigma_i$ can not intersect any reducing curve
in the Thurston canonical form for $\theta$.  It follows that there
are no reducing curves and the Thurston classification theorem
completes the proof.  \endproof

\section{The invariant $w(f,g).$}  

In this section we record some properties of a well known $\Z/2\Z$
invariant $\w(f,g)$ associated to a pair of commuting elements of $f,g
\in \Homeo_+(S^2)$.  Suppose $f,g \in \Homeo_+(S^2)$ commute and $f_t,
g_t$ are isotopies from $f = f_1$ and $g = g_1$ to $id = f_0 = g_0.$
We do not assume that $f_t$ and $g_t$ commute for $t \ne 0,1.$ The
invariant $w(f,g) \in \pi_1(\Homeo_+(S^2))$ is the
element determined by the loop $\gamma(t) =[f_t, g_t] =
f_t^{-1}g_t^{-1}f_t g_t$ in $\Homeo_+(S^2)$ with basepoint $id.$ We
recall that $\Homeo_+(S^2)$ is homotopy equivalent to $SO(3)$ \cite{Kn} and hence that that $\pi_1(\Homeo_+(S^2)) \cong \Z/2\Z$.

\begin{lemma}\label{lem:w(f,g)}
  The invariant $w(f,g) \in \pi_1(\Homeo_+(S^2))$ is also
represented by the loops $\alpha(t) =[f, g_t]$ 
 and $\beta(t) =[f_t, g]$ in $\Homeo_+(S^2)$ with basepoint $id.$
As a consequence it is independent of the choice of isotopies
$f_t$ and $g_t$ 
\end{lemma} 

\begin{proof}
Consider the function $\phi(s,t) =[f_s, g_t]$ defined on the square
$0 \le s \le 1,\ 0 \le t \le 1$ with values in $\Homeo_+(S^2)$.
The loop $\gamma(t) = \phi(t, t)$ is homotopic to 
\[
\alpha_0(t) = 
\begin{cases}
  \phi(2t,0) \text{ when } 0 \le t \le 1/2,\\
  \phi(1,2t-1) \text{ when } 1/2 \le t \le 1.\\
\end{cases}
\]
But $\alpha_0(t)$ is just a reparametrization of $\alpha$ which is
constant on the interval $0 \le t \le 1/2,$ so the loop $\gamma(t)$
is homotopic to $\alpha(t).$  The fact that $\gamma(t)$
is homotopic to $\beta(t)$ is proved similarly.
\end{proof}

\begin{prop}  Suppose $f,g,$ and $h$ are commuting elements 
of $\Homeo_+(S^2)$.  Then $w(fg, h) = w(f, h) + w(g, h)$ and
$w(f,g) = w(g,f).$
\end{prop} 
\begin{proof}
Consider the function $\phi(s,t) =[f_sg_t, h]$ defined on the square
$0 \le s \le 1,\ 0 \le t \le 1$ with values in $\Homeo_+(S^2)$.
The loop $\phi(t, t)$ represents $w(fg, h)$ by Lemma \ref{lem:w(f,g)}.
This loop is homotopic to the loop
\[
\alpha(t) = 
\begin{cases}
  \phi(2t,0) \text{ when } 0 \le t \le 1/2,\\
  \phi(1, 2t-1) \text{ when } 1/2 \le t \le 1.\\
\end{cases}
\]
But $\alpha(t)$ is the concatenation of the two loops
\[
\phi(2t,0) = g_0^{-1}f_{2t}^{-1} h^{-1} f_{2t} g_0 h
= f_{2t}^{-1} h^{-1} f_{2t} h = [f_{2t},h],
\]
where $0 \le t \le 1/2,$ and
\[
\phi(1, 2t-1 ) = g_{2t-1}^{-1}f_{1}^{-1} h^{-1} f_{1} g_{2t-1} h
= g_{2t-1}^{-1} h^{-1}  g_{2t-1}h = [g_{2t-1}, h],
\]
where $1/2 \le t \le 1.$ 
The first of these loops represents $w(f,h)$ and the second
represents $w(g,h)$. Similarly one can show $w(h,fg) = w(h,f) + w(h,g).$

To see the symmetry of $w(f,g)$ observe that for any $h$
we trivially have $w(h,h) = 0.$ Hence expanding 
$w(fg, fg)$ using the bilinearity just established we get
\[
0 = w(fg, fg) = w(f,f) + w(f,g) +w(g,f) + w(g,g) =  w(f,g) +w(g,f).
\]
Since $w$ takes values in $\Z/2\Z$ this implies $w(f,g) = w(g,f).$
\end{proof}

The universal covering space of $\Homeo_+(S^2)$ which we denote
$\widetilde \Homeo_+(S^2)$ is itself a group.  Since
$\pi_1(\Homeo_+(S^2)) \cong \Z/2\Z$ this is a two-fold covering.  

\begin{prop} If $f,g \in \Homeo_+(S^2)$ commute then
$w(f,g) = 0$ if and only if there are lifts
$\ti f, \ti g \in \widetilde \Homeo_+(S^2)$ of $f$ and
$g$ respectively such that $\ti f \ti g = \ti g \ti f.$ 
\end{prop}

\begin{proof}
Let $\ti f, \ti g \in \widetilde \Homeo_+(S^2)$ be lifts of $f$ and
$g$ respectively.  By definition $\ti f$ is a homotopy class of paths
from $id$ to $f$ in $\Homeo_+(S^2)$ relative to the endpoints.  We let
$f_t$ be a representative path for $\ti f$ and $g_t$ a representative
path for $\ti g.$ Then $w(f,g)$ is the homotopy class of the loop
$\alpha(t) = [f_t,g_t]$ in $\Homeo_+(S^2)$.  This loop is null
homotopic if and only if $\alpha$ lifts to a closed loop $\ti
\alpha(t)$ in $\widetilde \Homeo_+(S^2)$ with endpoints at
$\widetilde{id}$, the identity element of the group $\widetilde
\Homeo_+(S^2)$.  But this occurs if and only if $[\ti f, \ti g] =
\widetilde{id}$, i.e., if and only if $\ti f \ti g = \ti g \ti f.$
\end{proof}

We recall the following result from Lemma~1.2 of \cite{han:commuting}.

\begin{lemma} \label{standard example} 
If $f \in \Homeo_+(S^2)$ is conjugate to a non-trivial rotation and $g
\in \Homeo_+(S^2)$ commutes with $f$ and permutes the elements of
$\Fix(f)$ then $\w(f,g)=1$.
\end{lemma}

Suppose that $N$ is a  genus zero surface  and that $h : N \to N$ is a homeomorphism. Let $\bar N$ be the closed surface obtained from the interior of $N$ by compactifying each of its  ends with a point.  Then $\bar N$ is naturally identified with $S^2$ and  there is an induced homeomorphism $\bar h : S^2 \to S^2$.  We allow the possibility that $N = S^2$ in which case $\bar h= h$.

The following result is essentially  Lemma~3.1 of \cite{han:commuting}.  We include a proof for the reader's convenience.

\begin{lemma} \label{same for Thurston form} Assume the notation of Proposition~\ref{normal form} with  $S= \R^2$ or $S = S^2$.   Denote $\langle \theta_1,\ldots,\theta_n\rangle$ by $\Theta$.   Suppose that  $M_i$ has finite type and is $\Theta$-invariant, that $f \in \F$ corresponds to $\psi \in \Theta$ and that $g \in \F$ corresponds to $\phi \in \Theta$.  Then  $w(\bar f, \bar g) = w(\overline{\psi|_{M_i}},\overline{\phi|_{M_i}})$.   
\end{lemma}

\proof Denote $\psi|_{M_i}$ by $\psi_i$ and $\phi|_{M_i}$ by $\phi_i$.  Let $X$ be a finite set of punctures in $M$, one in each component of the complement of $M_i$, and let $\bar X$ be the corresponding set in $\bar S = S^2$.    Choose isotopies $\bar f_t$ and $\bar g_t$ rel $\bar X$, $t \in [1,2]$,   such that:
\begin{itemize}
\item [(1)]$\bar f_1 = \bar f$ and  $\bar f_2=\bar \psi$.
\item [(2)]$\bar g_1 = \bar g$ and  $\bar g_2=\bar \phi$.
\end{itemize}

  Let $M_i'$ be the complement in $M_i$ of a regular neighborhood of  $\partial M_i$.   Choose isotopies $\bar f_t$ and $\bar g_t$ rel $M_i'$, $t \in [2,3]$,   such that:
\begin{itemize}
\item [(3)]$\bar f_2 = \bar \psi$ and  $\bar f_3=\overline{\psi_i}$.
\item [(4)]$\bar g_2 = \bar \phi$ and  $\bar g_3=\overline{\phi_i}$.
\end{itemize}
Define $\bar h_t = [\bar f_t,\bar g_t]$ for $t \in [1,3]$.  By the third item of Proposition~\ref{normal form},  $\bar h_t$ is a loop in  $\Homeo_+(S^2)$ with basepoint $id$ and it suffices to show that this loop represents the trivial element of $\pi_1(\Homeo_+(S^2)$.  Moreover $\bar h_2|_{M_i}$ is the identity.    Choose an isotopy $\bar q_t$ rel $M_i' \cup X$, $t \in [0,1]$,  such that:
\begin{itemize}
\item [(5)] $\bar q_0 = \bar h_2$ and $\bar q_1 = $ identity.
\end{itemize}

Define loops in $\Homeo_+(S^2)$ by
\[
\bar \alpha(t) = 
\begin{cases}
  \bar h_t \text{ when } 1 \le t \le 2,\\
  \bar q_{t-2} \text{ when } 2 \le t \le 3\\
\end{cases}
\]
and 
\[
\bar \beta_t = 
\begin{cases}
  \bar q_{2-t} \text{ when } 1 \le t \le 2,\\
  \bar h_t \text{ when } 2 \le t \le 3.\\
\end{cases}
\]
Then $\bar h_t$ and the concatenation of $\bar \alpha_t$ and $\bar
\beta_t$ determine the same element of \newline
$\pi_1(\Homeo_+(S^2))$.  Since $\bar \alpha_t$ and $\bar \beta_t$ are
both relative to at least three points, Lemma~1.1 of
\cite{han:commuting} implies that both $\bar \alpha_t$ and $\bar
\beta_t$ determine the trivial element of $\pi_1(\Homeo_+(S^2))$.
\endproof

The following lemma will be used repeatedly in the proofs of Theorem~\ref{thm: partial plane} and   Theorem~\ref{thm: main sphere}.  

\begin{lemma}  \label{moving on} 
Assume that $M$ is a genus zero surface with finite type and that
$\Theta$ is a finite or virtually cyclic  abelian subgroup of $\Homeo_+(M)$,
each element of which either has finite order or is pseudo-Anosov.
Then there is a subgroup $\Theta'$ of $\Theta$ with index at most two
such that the sum of the number of $\Theta'$-fixed punctures, the
number of $\Theta'$-fixed points in the interior of $M$ and the number
of $\Theta'$-invariant components of $\partial M$ is at least
two. Moreover, if one of the following conditions is satisfied
\begin{itemize}
\item there is at least one $\Theta$-fixed puncture or
$\Theta$-invariant component of $\partial M$
 \item $\w(\phi, \psi)=0$ for all $\phi,\psi \in \Theta$.
\item $\Theta$ is  cyclic  
\end{itemize}
 then one may choose $\Theta' = \Theta$:
\end{lemma}

\proof Let $\{\theta_j\}$ be a set of generators for $\Theta$ and let
$\bar \Theta$ be the abelian subgroup of $\Homeo_+(\bar M=S^2)$ determined
by $\Theta$.  Suppose at first that some non-trivial $\bar \theta \in
\bar \Theta$ has finite order and so is conjugate to a rotation (see for example Remark~2.4 of \cite{han:commuting}).
Since $\bar \Theta$ is abelian, it acts on the two point set
$\Fix(\bar \theta)$.  Let $\bar \Theta'$ be the subgroup of $\bar
\Theta$ for which this action is trivial.  If the first bulleted item
is satisfied then $\bar \Theta' = \bar \Theta$ because a permutation
of a two point set that fixes at least one point fixes both
points. Similarly if the second bulleted item is satisfied then Lemma~\ref{standard example}
implies that $\bar \Theta' = \bar \Theta$. If the third bulleted item is satisfied then it is obvious that $\bar \Theta' = \bar \Theta$.  The lemma follows in this
case from the fact that a point in $\bar M$ corresponds to either a
puncture in $M$, a boundary component of $M$ of a point in the
interior of $M$.

   If every non-trivial element of $\bar \Theta$ has infinite order then $\bar \Theta$ is an infinite cyclic group (see, for example, Lemma 2.3 of \cite{han:fpt}) generated by a pseudo-Anosov homeomorphism $\bar \psi$.  The fixed points of $\bar \psi$ have index at most one so there are at least two of them and each is fixed by $\bar \theta$ for all $\theta$. The proof now concludes as in the previous case.
\endproof

\begin{lemma}\label{lem:gl_loop}
Let $\alpha_t$ be a loop in $\gl2$ with $\alpha_0 = \alpha_1 = I.$
For $y \in \R^2$ let $ev_y : \gl2 \to \R^2$ be the evaluation
map defined by $ev_y(A)= Ay.$
Suppose that there is $x  \in \R^2 \setminus\{0\}$ such that
the loop $ev_{x}(\alpha_t) = \alpha_t x$ is 
inessential in $\R^2 \setminus\{0\}$.
Then $\alpha_t$ is inessential in $\gl2$.
\end{lemma}

\begin{proof}
There is a deformation retraction $\cR$ of $\gl2$ onto $\so2,$
the group of rotations of the plane.  This follows from the matrix
polar decomposition according to which any real matrix $A$ can be
written uniquely as $PO$ with $P$ symmetric positive definite and $O$
orthogonal (see \cite{W} pp. 131-6).  Since this decomposition is
continuous and the symmetric positive definite matrices are convex
this gives $\cR.$

Let $\beta_t = \cR \circ \alpha_t$ so $\beta_t$ is
loop in  $\so2 \subset \gl2$  and the loop
$ev_{x}(\beta_t) = \beta_t x$ is 
is inessential.
The evaluation map $ev_{x} : \so2 \to \R^2$
is a homeomorphism from $\so2$ onto the circle of radius
$r = \| x \|$ in $\R^2$.  It follows that $\beta_t$ is an
inessential loop in $\so2$ and hence that $\alpha_t$ is
inessential in $\gl2.$
\end{proof}

\begin{prop}
Suppose $f$ and $g$ are commuting elements of $\Diff^1_+(\R^2)$
and $\bar f$ and $\bar g$ are their extensions to $S^2,$ the one-point
compactification of $\R^2.$  Then $w(\bar f,\bar g) = 0.$
\end{prop}
\begin{proof}
Clearly it suffices to prove that the loop $h_t = [f_t, g_t]$
in $\Diff^1_+(\R^2)$ is contractible.  

Choose balls $B_1 \subset B_2$ in $\R^2$ centered at the origin
and a function $H: [0,1]^2 \to \Diff^1_+(\R^2)$ such that
\begin{enumerate}
\item $h_t(0) \in B_1$ for all $0 \le t \le 1.$
\item $H(s,t)(x) = x - sh_t(0)$ for all 
$x \in B_1$ and all $0 \le s,t \le 1$, and 
\item $H(s,t)(x) = x$ for all 
$x$ outside $B_2$ and all $0 \le s,t \le 1$.

\end{enumerate}
Then if we define $\theta_t = H(1,t) \circ h_t$, the loop $\theta_t$
in $\Diff^1_+(\R^2)$ is homotopic to $h_t$ (the homotopy being given
by $H(s,t) \circ h_t,\ 0\le s \le 1$.  Moreover, $\theta_t(0) = 0$ for all $t
\in [0,1]$ and $\theta_t(x) = h_t(x)$ for $x$ outside $B_2.$
To prove our result it suffices to show $\theta_t$ is an
inessential loop in $\Diff^1_+(\R^2).$

To do this we will apply 
a result of D. Calegari (see Section 4.1 and
Theorem D from Section 4.8 of \cite{Calegari}).  This result asserts that
if $G$ denotes the group generated by $f$ and $G_\infty$ is the
group of germs at $\infty$ of $G$ then there is a split exact
sequence
\[
0 \to \Z \to \hat G_\infty \to G_\infty \to 0
\]
where $\hat G_\infty$ is the group of all lifts of $f$ and $g$
to the universal cover $\ti A$ of the germ $A$ of a neighborhood 
of $\infty.$ The fact that this sequence splits implies there
are lifts $\ti f, \ti g \in \ti G_\infty$ of the germs of $f$ and $g$ 
respectively such that $[\ti f, \ti g] = id.$ 

This implies that if $N$ is a sufficiently small punctured disk
neighborhood of $\infty$ and $x_0 \in N$ then 
$h_t(x_0)$ is an inessential loop in $N.$   Moreover
if $N$ is sufficiently small then 
$\theta_t(x_0)$ is an inessential loop in $N$ since it
agrees with $h_t(x_0)$.

We now perform a further homotopy 
from the loop $\theta_t$ in $\Diff^1_+(\R^2)$ to 
the loop $\phi_t = \Phi(0,t)$ where 
$\Phi(s,t) \in \Diff^1_+(\R^2)$  is defined by
\[
\Phi(s,t)(x) =
\begin{cases}
\frac{\theta_t(sx)}{s}  \text{ when } 0 < s \le 1,\text{ and}\\
D\theta_t(0)(x) \text{ when } s =0,\\
\end{cases}
\]

It now suffices to show the loop $\phi_t$ is inessential in 
$\Diff^1_+(\R^2).$  We observe that the loop
$\phi_t(x_0)$ is homotopic to $\theta_t(x_0)$ in
$\R^2\setminus\{0\}$ and therefore inessential.  We also
note that $\phi_t = D\theta_t(0) \in \gl2$.  Hence
Lemma\ref{lem:gl_loop} implies the loop $\phi_t$ is inessential in 
$\Diff^1_+(\R^2).$
\end{proof}

As an immediate corollary we have the following.
\begin{cor}\label{cor: fp=>w=0}
Suppose $f$ and $g$ are commuting elements of $\Diff^1_+(S^2)$
which have a common fixed point $x_0$.  Then 
$w(f,g) = 0.$
\end{cor}

\section{Nielsen Classes}\label{sec:nielsen}
We recall some definitions from Nielsen fixed point theory that are
needed in the proof of Proposition~\ref{nonzero nielsen class}.
Further details can be found for example in Brown \cite{Brown}.

Let $h : N \to N$ be a homeomorphism of a connected surface $N$.  In the paragraphs preceding Lemma~\ref{max exists} we recalled the definition of Nielsen equivalence and  the relationship between Nielsen classes and and lifts of $h$ to the universal cover of $N$.

We say that a Nielsen class $\mu$ for $h : N \to N$ is {\em compactly covered} if $\Fix(\ti h)$ is
compact for some, and hence every, lift $\ti h$ for $\mu$.  In this
case, $\Fix(\ti h)$ is contained in the interior of a compact disk $D
\subset \ti N$ and the {\em Nielsen index $I(h,\mu)$} for $\mu$ is the
usual index obtained as the winding number of the vector field
$x - \ti h(x)$ around $\partial D$ (see, e.g. Brown \cite{Brown}).

If $Q$ is an isolated puncture in $N$ then we say that a  Nielsen class $\mu$ in $\Fix(h)$ {\em peripherally contains} $Q$ if for some, and hence every,  $x$ in $\mu$  there is a properly immersed ray $\alpha \subset N$ based at $x$ that converges to $Q$ such that   $h(\alpha)$ is properly homotopic to $\alpha$ relative to $x$. 

  Given an isotopy $h_t$  from $h_0 = h$ to some $h_1$, let $\ti h$ be a lift for $\mu$ and let $\ti h_t$ be the lift of $h_t$ that begins with $\ti h_0 = \ti h$.  The  Nielsen class $\mu_t$  for $h_t$ determined by $\ti h_t$ is independent of the choice of $\ti h$.  We say that {\em $\mu_t$ is the Nielsen class determined by $\mu$ and $h_t$}.

\begin{lemma}  \label{nielsen}Suppose that $\mu$ is a non-trivial, compactly covered Nielsen class for $h : N \to N$ and that each fixed end of $N$ is isolated and not peripherally contained in $\mu$.  If $h_t$ is an isotopy with $h_0 = h$ and $\mu_t$ is determined by $h_t$ and $\mu = \mu_0$, then    each $\mu_t$  is compactly covered, does not peripherally contain any fixed end  and $I(h_t,\mu_t)$ is independent of $t$.
\end{lemma}

\proof Choose $\ti h_t$ as in the definition of $\mu_t$.  To prove
that $\mu_t$ is compactly covered and that $I(h_t,\mu_t)$ is
independent of $t$, it suffices to show that $\ti X := \cup_t \Fix(\ti
h_t)$ is compact in $\ti N$.  One can then compute $I(h_t,\mu_t)$
using a disk $D \subset \ti N$ whose interior contains $\Fix(\ti h_t)$
for all $t$ and it is then a standard fact (see Dold \cite{Dold})
that $I(h_t,\mu_t)$ is independent of $t$.

Let $Q_1,\ldots,Q_k$ be the isolated fixed ends of $N$.  Choose neighborhoods $V_i \subset U_i \subset N$ of $Q_i$ such that $h_t(V_i) \subset U_i$ for all $t$ and such that the one point compactification of $U_i$ obtained by filling in the puncture $Q_i$ is a disk.  Then $Q_i$ is peripherally contained in any   Nielsen class for $h_t$ that contains an element of $\Fix(h_t) \cap V_i$.      Choose a neighborhood $W$ of the remaining ends of $N$ so that $\Fix(h_t) \cap W = \emptyset$ for all $t$.

Define $H : N \times [0,1] \to N \times [0,1]$ by $H(x,t) =
(h_t(x),t)$ and $\ti H : \ti N \times [0,1] \to \ti N \times [0,1]$ by
$\ti H(x,t) = (\ti h_t(x),t)$.  Let $Y \subset \Fix(H)$ be the
projected image of $\ti Y:= \Fix(\ti H)$.  Any two points $(x_s,s),
(x_t,t) \in Y$ are connected by an arc $\gamma \subset N \times [0,1]$
such that $H(\gamma)$ is homotopic to $\gamma$ rel endpoints.  Since
$\mu_0 = \mu$ is non-trivial and compactly covered, $\ti h_0 = \ti h$
does not commute with any covering translations. It follows that $\ti
H$ does not commute with any covering translations and hence that each
element of $Y$ has a unique lift to an element of $\ti Y$.  We show
below that $Y$ is compact.  Since points in $Y$ lift uniquely there is
a continuous lifting map from $Y$ to $\ti Y$.  It follows that $\ti
Y$, and hence $\ti X$, is compact.

Denote the set of elements of the equivalence class $\mu_t$ by $X_t$.  Then $X :=\cup_tX_t \subset N$  is the projected image of $\ti X \subset \ti N$ and of $Y \subset N\times [0,1]$.  
 Obviously $X \cap W = \emptyset$.  To complete the proof that $Y$ is compact, it remains to check that  $X_t \cap V_i = \emptyset$ for all $i$ and $t$.   We will prove the stronger statement that $\mu_t$ does not peripherally contain any puncture, which  completes the proof of the lemma.  Suppose to the contrary that there is a ray $\alpha_t \subset N$  connecting $x_t$ to $Q_i$ such that $h_t(\alpha_t)$ is properly homotopic to $\alpha_t$ rel $x_t$.   For  any $x_0 \in X_0$, there is an arc $ \gamma \subset  N \times [0,1]$ connecting $(x_0,0)$ to $(x_t,t)$ such that $H( \gamma)$ is homotopic to $ \gamma$ rel endpoints.   Let $\alpha_0 \subset N$ be the concatenation of the projection $\gamma_0$ of  $ \gamma$ and  $\alpha_t $.   Then $h_0(\alpha_0)$ is properly homotopic to $\alpha_0$ rel $x_0$.  This implies that $\mu_0$ peripherally contains $Q_i$ which is a contradiction.    
\endproof

Suppose that $\theta :S\to S$ is an orientation preserving
homeomorphism of a finitely punctured surface, that $C$ is a compact
set in the complement of $\Fix(\theta)$ and that $x \in \Fix(\theta)$.
Let $N$ be the component of $S \setminus C$ that contains $x$ and let
$h = \theta|_N : N \to N$.  Define {\em the Nielsen class for $\theta$
relative to $C$ determined by $x$} to be the Nielsen class for $h$
determined by $x$. The remaining \lq relative\rq\ definitions are made
similarly, using $h$ in place of $\theta$.

The following proposition  is the main result of this section.   

\begin{proposition} \label{nonzero nielsen class}  Suppose that $\theta: S \to S$ is an orientation preserving homeomorphism of  a finitely punctured surface, that  $C$ is a compact invariant set in the complement of $\Fix(\theta)$, that $M_i$ is an essential  subsurface of $M :=S \setminus C$  with finite negative Euler characteristic and that $\theta|_{M_i}$ is either pseudo-Anosov relative to its puncture set or has finite order greater than one.  Suppose further that $z \in \Fix(\theta)$ is  contained in the interior of $M_i$ and that     $f_t :S \to S$ is an  isotopy rel $C$ from  $f_0 =\theta$ to some   $f =f_1$.   Let $\mu$ be the Nielsen class for $\theta$ relative to $C$ determined by $z$ and let $\nu$ be the Nielsen class for $f$ relative to $C$ determined by $\mu$ and $f_t$.  Then $\nu$ is non-trivial, compactly covered  and does not peripherally contain any punctures.  
\end{proposition}

\proof  There is no loss in assuming that $M$ is connected.  Choose   lifts $\ti z$ and $\ti \theta$ such that $\ti z \in \Fix(\ti \theta)$.   Since $M_i$ is essential, the component of the pre-image of $M_i$ in $\ti M$ that contains $\ti z$ is a copy of the universal cover $\ti M_i$ of $M_i$.   The following properties of $\theta|_{M_i}$ and $\ti \theta|_{\ti M_i}$ are well known.
\begin{itemize}
\item  $\Fix(\ti \theta|_{\ti M_i}) = \{\ti z\}$ and the  fixed point index of $z$ with respect to $\theta$ is non-zero.
\item   No component of $\partial \ti M_i$, and hence no component of $\ti N \setminus \ti M_i$, is $\ti \theta$-invariant.
\item  There does not exist a ray $\alpha$ in $M_i$ connecting $z$ to a puncture of $M_i$ such that $\theta(\alpha)$ is homotopic to $\alpha$ rel endpoints.
\end{itemize}
The first two items imply that $\Fix(\ti \theta) = \{\ti z\}$ and hence that $\mu$ is compactly covered. Suppose that there is a  ray $\alpha$ in $M$ connecting $z$ to a puncture of $M$ such that $\theta(\alpha)$ is homotopic to $\alpha$ rel endpoints.  Decompose $\alpha$ into an alternating concatenation of maximal subpaths $\sigma_j$ in $M_i$ and maximal subpaths $\tau_j$ in the complement of $M_i$.   With the exception of $\sigma_1$ and the final $\tau_j$ if it is a ray, all of these subpaths have both endpoints in $\partial M_i$ and we may assume that they cannot be homotoped rel endpoints into $\partial M_i$.  In this case, the component of $\partial \ti M_i$ that contains the terminal endpoint of the lift $\ti \sigma_1$ that begins at $\ti z$ is preserved by $\theta$ which contradicts the second item above.  We conclude that $\mu$ does not peripherally contain any punctures.  The proposition now follows from Lemma~\ref{nielsen} and the first item above.
\endproof

\section{A Fixed Point Lemma}

   The Brouwer translation theorem implies that if $f \in \Homeo_+(\R^2)$ has a compact invariant set $C$ then $f$ has a fixed point.  We will enhance this by adding some more information about the location of the fixed point relative to $C$.

\begin{defn} \label{pivot} 
Suppose that $f \in \Homeo_+(\R^2)$ and that $C \subset \R^2$ 
is a compact invariant set with $C \cap \Fix(f) = \emptyset.$
 Define $P(C,f)$ to be the union of Nielsen classes of  $\Fix(f)$ rel $C$
which do not peripherally contain $\infty$.  
\end{defn}

We consider first the special case that $C$ is a periodic orbit.  The following
result is due to Gambaudo \cite{Gam}.  We provide a proof as a convenience
to the reader and because it is quite short given our earlier lemmas.

\begin{lemma}[Gambaudo]  \label{lem: periodic} 
If $f \in \Homeo_+(\R^2)$ has a periodic point $p$
of period greater than one, then $P(\Orb(p),f) \ne  \emptyset$.  
\end{lemma}

\proof Let $ N = \R^2 \setminus \Orb(p)$.  By the Thurston
classification theorem, there is a homeomorphism $\theta : N \to N$
that is isotopic to $f$ and an essential $\theta$-invariant subsurface
$N_1$ of negative Euler characteristic that contains $\infty$ as a
puncture such that $\theta|_{N_1}$ is either periodic or pseudo-Anosov.
No puncture of $N$ other than $\infty$ is fixed.  If $\partial N_1$ is
non-empty then it contains at least two components, each of which
bounds a disk in the complement of $N_1$ that contain punctures
corresponding to $\Orb(p)$.  Since the action of $\theta$ on $\Orb(p)$
is minimal, these disks are not $\theta$-invariant and neither are the
components of $\partial N_1$.  The lemma now follows from
Lemma~\ref{moving on} and Lemma~\ref{nonzero nielsen class}.
\endproof

We can now state and prove our general result on $P(C,f)$.

\begin{prop}  \label{prop: pivot} 
Suppose that $f \in \Homeo_+(\R^2)$ and that $C \subset \R^2$ 
is a non-empty compact $f$-invariant set with $C \cap \Fix(f) = \emptyset.$
Then $P(C,f)$ is non-empty and compact.  
\end{prop} 
\begin{proof}
The proof is by contradiction.
Suppose that  that each Nielsen class in$\Fix(f)$  peripherally contains 
  $\infty.$  Choose open neighborhoods $V_{\infty} \subset U_{\infty}$ of $\infty$ that are disjoint from $C$ and that satisfy  $f(V_{\infty}) \subset U_{\infty}$.      Similarly, for each $x \in \Fix( f)$ choose open disk neighborhoods $V_x \subset U_x$ whose closures are disjoint from $C$  with the property that 
$f(V_x) \subset U_x$.  Let $V_{\infty} \cup \{V_{x_i}\}$ be a finite subcover of the cover
$V_{\infty} \cup \{V_{x}\}$ of $\Fix(f)$.  For each $x_i$ there is a 
path from $\infty$ to $x_i$ which is homotopic rel endpoints to its $f$
image and the support of the homotopy is disjoint from $C$. For
each $x \in V_{x_i}$ we can form a path from $\infty$ to $x$
by concatenating the path from $\infty$ to $x_i$ with a path
in $V_{x_i}$ from $x_i$ to $x$.  From this it is clear that
for some $\epsilon >0$ and for any $x \in \Fix(f)$ 
there is a path $\alpha$ from $\infty$ to $x$ 
and a homotopy rel endpoints from $\bar f(\alpha)$ to $\alpha$ whose support
is disjoint from an $\epsilon$ neighborhood of $C$.  

Now let $z \in C$ be a recurrent point.  By a small $C^0$ perturbation
of $f$ supported in the $\epsilon$ neighborhood of $C$ we can
make $z$ be a periodic point.  Moreover, we may do this with a
perturbation supported on an open set $W$ disjoint from its $f$
image.  Hence the perturbed map has precisely the same fixed point set
as $f$ and every fixed point   still belongs to a   Nielsen class that peripherally contains
$\infty$ since $f$ was not modified outside an $\epsilon$
neighborhood of $C$.  This contradicts Lemma \ref{lem: periodic}.

This proves $P(C,f)$ is not empty.  It is immediate that $P(C,f)$ is
compact since any point of $\Fix(f) \cap V_{\infty}$ determines a Nielsen class that peripherally contains $\infty$.  
\end{proof}

\begin{remark} \label{pivot set is natural} It is immediate that if
$g \in \Homeo_+(\R^2)$ commutes with $f$ and $C$ is $g$-invariant
then $P(C,f)$ is also $g$-invariant.
\end{remark}

\section{Proof of Theorem~\ref{thm: partial plane}}

We will prove the following strengthened version of Theorem~\ref{thm: partial plane}.  The extra condition (2)  plays an essential role in our inductive proof.

\begin{thm}  \label{thm: main plane} 
Suppose that $\F$ is a finitely generated abelian
subgroup of $\Diff^1_+(\R^2)$  and that there is a compact $\F$-invariant set $C\subset \R^2$.  Then the following hold:
\begin{enumerate}
\item $\Fix(\F)$ is non-empty.
\item If $\F = \langle f_1, \dots, f_n \rangle$ and
$C \subset \Fix(f_1, \dots, f_{n-1})$ then
\[
(C \cup P(C, f_{n})) \cap \Fix(\F) \ne \emptyset.
\]
\end{enumerate}
\end{thm}

\proof
The proof is by induction on the rank $n$ of $\F$.   
Our strategy is  to
use the following intermediate inductive statement.
\begin{itemize}
\item[$A_n:$]   (1) holds for $\F = \langle f_1, \dots, f_n \rangle$  
under the additional hypothesis that $\Fix(f_n)$ is compact.
\end{itemize}

We also make use of the following notation.

\begin{itemize}
\item[$B_n:$] Theorem~\ref{thm: main plane} holds for $\F = \langle f_1, \dots, f_n \rangle$.
\end{itemize}

 Obviously,  the validity of $B_n$ for all $n>0$ is the result we wish to
prove.  $A_1$  follows from the Brouwer translation theorem 
and  $B_1$ is a consequence of Proposition~\ref{prop: pivot} (recall
that $\Fix(\emptyset) = \R^2$ so in $B_1$ the condition
that $C \subset \Fix(f_1, \dots f_{n-1})$ is trivially true).
Hence by induction the theorem follows from Lemmas \ref{lem: B=>A} 
and \ref{lem: AB=>B} below.

\begin{lemma}  \label{lem: B=>A} 
$B_{n-1} \Rightarrow A_{n}$ for all $n \ge 2$.
\end{lemma}

\proof  We  assume $B_{n-1}$, the hypothesis of $A_n$ and that  $\Fix(\F) = \emptyset$ and
show this leads to a contradiction.  Denote $\langle f_1, \dots, f_{n-1} \rangle$ by $\F_{n-1}$.  

The set  $L := \Fix(f_{n})$   is non-empty by
the Brouwer translation theorem and is compact by the  hypothesis of $A_{n}$.   
Define  
 $K = P(L, f_{n-1}) \cap \Fix(\F_{n-1})$.  Since $K$ is a subset of $P(L, f_{n-1})$ it is
compact. We claim that 
it is also non-empty.  To see this, observe that $J =
\Fix(f_1,\dots,f_{n-2},f_{n}) \ne \emptyset$ by  $B_{n-1}$  and is compact since it is a subset of $L$.
$B_{n-1}$ applied to $\F_{n-1}$ implies
\[
(J \cup P(J, f_{n-1})) \cap \Fix(\F_{n-1})  \ne \emptyset.
\]
If $J \cap \Fix(\F_{n-1})  \ne \emptyset$ then 
$\Fix(\F)  \ne \emptyset$, a contradiction. 
It follows, since $J \subset L,$  that 
\[
K = P(L, f_{n-1}) \cap \Fix(\F_{n-1})
 \supset P(J, f_{n-1}) \cap \Fix(\F_{n-1})  \ne \emptyset
\]
as desired.

Apply Proposition~\ref{normal form} to $\F, K$ and $L$ and assume the notation of that proposition.  Denote $\langle \theta_1,\ldots,\theta_n\rangle$ by $\Theta$ and $\langle \theta_1,\ldots,\theta_{n-1}\rangle$ by $\Theta_{n-1}$.  

We claim that the subsurface $M_1$ containing $\infty$ as a puncture
has at least one $\Theta$-invariant boundary component $\Gamma_1$.  If
either $\theta_n|_{M_1}$ is the identity or $\Theta_{n-1}|_{M_1}$ is
the identity then this follows from the fact that $K$ and $L$ are
disjoint and non-empty. Otherwise, $M_1$ is finitely punctured and the
claim follows from Lemma~\ref{moving on}, the assumption that
$\Fix(\F) = \emptyset$ and the fact that $\Fix(\theta_n) \cap
\Int(M_1) = \Fix(f_n) \cap \Int(M_1) = \emptyset$.

     $\Gamma_1$ is a boundary component of some annulus in $\A$.  The other boundary component of this annulus is a boundary component $\Gamma_2$ of some subsurface, say $S_2$, that is $\Theta$-invariant.  If $\partial M_2$ has a second $\Theta$-invariant component $\Gamma_3$ then one can repeat this operation to find a third $\Theta$-invariant subsurface and so on.  This process must eventually end so there exists a $\Theta$-invariant subsurface $M_i$ with exactly one $\Theta$-invariant boundary component $\Gamma$; moreover, the disk $D$ bounded by $\Gamma$ and containing $M_i$ does not contain $\infty$.

  If both $K \cap D \ne \emptyset$ and $L \cap D \ne \emptyset$ then the argument applied above to $M_1$ would apply to $M_i$ to produce a $\Theta$-invariant element of $R$ that we have assumed does not exist.   We may therefore assume that $K \cap D = \emptyset$ or $L \cap D = \emptyset$.    
If  $L \cap D = \emptyset$ then $K \cap D\ne \emptyset$ and every Nielsen class in $L$ rel $K \cap D$ peripherally contains $\infty$.  This contradicts  Proposition~\ref{prop: pivot}.

It remains to consider the case that $K \cap D = \emptyset$ and $L \cap D \ne \emptyset$.    We show below that $C' = L\cap D \cap \Fix(f_1, \dots, f_{n-2}, f_{n})$ is non-empty. Assuming this for the moment, we argue to a contradiction thereby finishing the proof of the lemma.  Apply  $B_{n-1}$ to $\F_{n-1}$ and $C'$  to conclude that $(C' \cup
P(C', f_{n-1})) \cap \Fix(\F_{n-1}) \ne \emptyset$ and hence, since $\Fix(\F) = \emptyset$, that   $P(C', f_{n-1}) \cap \Fix(\F_{n-1}) \ne \emptyset$.   This contradicts the fact that $P(C', f_{n-1}) \subset D \cap P(L,f_{n-1})$ and the assumption that  $K \cap D = \emptyset$.

We now  show that $C' \ne \emptyset$
 as promised.  Let $M'$ be the surface obtained from $M = \R^2 \setminus (K \cup L)$ by filling in the punctures corresponding to the compact set $L_D:=L \cap D$.   Thus $D \subset M'$ and   $L_D$ is $\F$-invariant. Let $\ti M'$ be the universal covering
space of $M'$, let $\ti D \subset \ti M'$ be a lift of $D$ and let $\ti L_D$ be the lift of $L_D$ that is contained in $\ti D$.  Note that $\ti L_D$ is homeomorphic to $L_D$ and is therefore compact.   For each $f \in \F$ there is a unique lift $\ti f : \ti M' \to \ti M'$ that setwise fixes $\ti L_D$.    The uniqueness of $\ti f$ implies that $\langle \ti f_1, \ldots,\ti f_{n-2},\ti f_n\rangle$ is abelian.  Applying $A_{n-1}$ (which follows from $B_{n-1}$) to this group and $\ti L_D$, we conclude that $\Fix(\ti f_1, \ldots,\ti f_{n-2},\ti f_n) \ne \emptyset$.  This set evidently projects into $\Fix(f_1, \dots, f_{n-2}, f_{n})$ and it projects into $\Fix(f_n|_{M'}) \subset L \cap D$.  It therefore projects into $C'$. 
\endproof

\begin{lemma}  \label{lem: AB=>B} 
$A_{n}$ together with $B_j$ for $j \le n-1$ 
implies $B_{n}$ for all $n \ge 2$.
\end{lemma}
\proof   Let us first reduce to the case that $C \subset \Fix(\F_{n-1})$ by producing   a non-trivial compact $f_n$-invariant subset $C_n$ of $\Fix(\F_{n-1})$.  If  $C$ is not a subset of $\Fix(\F_{n-1})$ then replace it with  $C_n$.  

We prove by induction on $2 \le k \le n$ that there is a non-trivial
compact $\F_k$-invariant subset $C_{k}$ of $\Fix(
f_1,\ldots,f_{k-1})$. For $k = 2$, this follows from Proposition
\ref{prop: pivot} applied to $C$ and Remark~\ref{pivot set is
natural}.  Assume now that $C_{k-1}$ exists.  If $C_{k-1} \cap
\Fix(f_k) \ne \emptyset$ define $C_k = C_{k-1} \cap \Fix(f_k)$;
otherwise define $C_{k} = P(C_{k-1}, f_{k-1}) \cap \Fix(
f_1,\ldots,f_{k-1})$.  Then $C_k$ is non-empty by $B_{k-1}$, is
compact because $P(C_{k-1}, f_{k-1}))$ is compact and is
$\F_k$-invariant by Remark~\ref{pivot set is natural}.  This completes
the induction step and the proof that $C_n$ exists.
Hence we may now assume that $C \subset \Fix(\F_{n-1})$.  
 
We will establish $B_{n}$ if we show $(C \cup P(C, f_{n})) \cap
\Fix(\F_n) \ne \emptyset.$ If $C \cap \Fix(f_{n}) \ne \emptyset$ we
are done so we may assume $C$ is disjoint from $\Fix(f_{n}).$

Let $K$ be a subset of $C$ which is $f_{n}$ minimal.  Assume the
notation of Proposition~\ref{normal form} with $L = \emptyset$ and
$M = \R^2\setminus K.$ Denote
$\langle \theta_1,\ldots,\theta_n\rangle$ by $\Theta$ and let $M_1$
denote the subsurface which contains $\infty$ as a puncture.  There
are two possibilities.

\noindent
{\bf Case 1. }{\em $M_1$ has finite type:} \\ If $M_1$ has punctures
corresponding to elements of $K$ then $K$ is finite, $M = M_1$, and
the set of reducing curves $R = \emptyset$.  Otherwise, $\infty$ is
the only puncture in $M_1$ and the components of $\partial M_1$, of
which there are at least two, are transitively permuted by
$\theta_n$. In either case, $\theta_n$ is not the identity and
Lemma~\ref{moving on} implies that there is a $\Theta$-fixed point $x$
in the interior of $M_1$.  Let $\ti M$ be the universal cover of $M$,
let $\ti x \in \ti M$ be a lift of $x$, let $\ti M_1$ be the component
of the full pre-image of $M_1$ that contains $\ti x$ and let $\ti f_j$
be the lift of $f_j$ that is equivariantly isotopic to the lift $\ti
\theta_j$ of $\theta_j$ that fixes $\ti x$.  Since $\F = \langle f_j
\rangle$ and $\langle \theta_j|_{M_1}\rangle$ are abelian, so is $\ti
\F = \langle \ti f_j\rangle$.

 The lifts $\ti \theta_j$  of $\theta_j$ that fix $\ti x$  generate   an abelian subgroup  $\ti \Theta$  of $\Diff^1_+(\ti M)$.  Each $f_j$ has a unique lift $\ti f_j$ that is equivariantly isotopic to $\ti \theta_j$ and  $\ti \F := \langle \ti f_1,\ldots, \ti f_n\rangle$ is an abelian subgroup of $\Diff^1_+(\ti M)$.  

  By Proposition~\ref{normal form}-(3), $\theta_n|_{M_1}$ is either finite order or pseudo-Anosov.  Let $\mu$ be the $f_n$-Nielsen class relative to $K$ determined by $\ti f_n$.  Proposition~\ref{nonzero nielsen class}   implies that $\mu$ is non-empty, compact and does not peripherally contain $\infty$.  In particular, all elements of $\mu$ are contained in $P(K, f_{n}) \subset P(C, f_{n})$.  
   Since $\ti \F_n$ is abelian, $\Fix(\ti f_n)$ is $\ti \F$-invariant.  $A_{n}$ therefore implies that 
  some element of $\mu$ is contained in $\Fix(\ti \F).$  
 This verifies $B_n$ in case 1.

{\bf Case 2. }{\em $M_1$ is infinitely punctured:} \\
n this case, every element of $\Theta_{n-1}$ acts as the identity on
$M_1$.  Since $K$ is $\theta_n$-minimal, $M= M_1$ and $R = \emptyset$.
If there is an essential finite type subsurface $N \subset M_1$ that
is preserved up to isotopy by $\theta_n$, then, after modifying
$\theta_n$ by an isotopy, we can redefine $R$ to be $\partial N$. This
redefines $M_1$, and as just observed, the new $M_1$ must be finitely
punctured.  Case 1 completes the proof in this case so we may assume
that no such $N$ exists.

   Let $\ti M$ be the universal cover of $M = M_1$, let $\ti x \in \ti M$ be a lift of  $x \in P(K, f_{n})$ and  let $\ti f_{n}: \ti M \to \ti M$
be the lift of $f_{n}$ which fixes $\ti x.$    Lemma~\ref{augmented regular nbhd} implies that   $f_n$ does not fix the homotopy class of any essential non-peripheral closed curve.  It follows that  $\ti f_n$ does not commute with any covering translations.  Thus $\Fix(\ti f_n)$ is compact.

Since $\Theta_{n-1}$ is the trivial group, each $f \in \F_{n-1}$ is isotopic to the identity rel $K$.  There is a unique lift $\ti f : \ti M \to \ti M$ that is equivariantly isotopic to the identity.  The uniqueness of $\ti f$ implies that $\ti F_{n} = \langle\ti f_1,\ldots,\ti f_{n}\rangle$ is   abelian.   $A_{n}$ therefore implies that 
  there is a point $\ti z$ in $\Fix(\ti \F).$  
The projection of this point into $M$ is a point
$z \in \Fix(\F)$ in the same  Nielsen class relative to $f_n$ as $x$.  Thus $z \in P(K, f_{n}) \subset P(C, f_{n})$.  This verifies $B_n$ in case 2.
\endproof

\section{Proof of Theorem~\ref{thm: main sphere}}

We begin with a pair of applications of Theorem~\ref{thm: partial plane}, the first of which requires some technical notation.

Suppose that $S$ is a finitely punctured surface, that $\F$ is a
finitely generated abelian subgroup of $\Homeo_+(S)$ and that $C
\subset S$ is compact and $\F$-invariant.  Suppose further that  $M_1$ is an essential
subsurface of $M :=S \setminus C$ that has finite negative Euler
characteristic and that is $\F$-invariant up to isotopy.  There is a
homomorphism $\Phi$ from $\F$ to the mapping class group of $M_1$
defined by $f \mapsto [\theta|_{M_1}]$ where $\theta \in \Homeo_+(S)$
is isotopic to $f$ and preserves $M_1$ and where $[h]$ is the isotopy
class of $h$.  This is well defined because homeomorphisms of $M_1$
that are isotopic in $M$ are isotopic in $M_1$.       If each element of the abelian subgroup 
$\Theta_1^*:=\Phi(\F)$ is either of pseudo-Anosov type or has finite order
then we say that {\em $\F$ is $M_1$-irreducible}.       For example, if some element of 
$\Theta_1^*$ has pseudo-Anosov type then
Lemma~\ref{virt_cyclic} implies that $\F$ is $M_1$-irreducible.

      If $\F$ is $M_1$-irreducible then
Lemma~\ref{virt_cyclic} implies that $\Theta^*_1$ is
either finite or virtually cyclic and lifts to a subgroup $\Theta_1 \subset \Homeo_+(M_1)$.  By
construction, each $f \in \F$ is isotopic to some $\theta \in
\Homeo_+(S)$ where $M_1$ is $\theta$-invariant and $\theta|_{M_1} \in
\Theta_1$.

\begin{proposition} \label{interior fixed point}  
Assume notation as above and further that $\F$ is $M_1$-irreducible,
that $\Theta_1^*$ is non-trivial and that $\Theta_1 \subset \Homeo_+(M_1)$
is a lift of $\Theta_1^*$.  For each  $z \in \Fix(\Theta_1)$ that is contained in
the interior of $M_1$,   there exists $x \in \Fix(\F) \cap M$ with
the following properties: if $\theta_1^*:=\Phi(f)$ is non-trivial then
the $f$-Nielsen class rel $C$ that contains $x$ is compactly covered,
does not peripherally contain any punctures and corresponds, via the
isotopy between $f$ and $\theta$, to the $\theta$-Nielsen class rel
$C$ that contains $z$, where $\theta \in
\Homeo_+(S)$ is isotopic to $f$  and $\theta|_{M_1} \in
\Theta_1$.
\end{proposition}

\proof We may assume without loss that $M$ is connected.  Let $\ti
M_1$ be the universal cover of $M_1$ and let $\ti z \in \ti M_1$ be a
lift of $z$.  Each element of $\Theta_1$ has a unique lift to $\ti
M_1$ that fixes $\ti z$.  This lifts the group $\Theta_1$
to a group $\ti \Theta_1 \subset \Homeo_+(\ti M_1).$

Lift $\F$ to $\ti \F \subset \Homeo_+(\ti M)$ as follows.  
Given $f \in \F$, let $\ti f_1 : \ti M_1 \to \ti M$
be the unique lift of $f|_{M_1}$ that is equivariantly isotopic to an
element of $\ti \Theta_1$ and let $\ti f : \ti M \to \ti M$ be the
unique extension of $\ti f_1$ to a lift of $f$.  The uniqueness
properties imply that $f \mapsto \ti f$ defines a lift of the
group $\F.$

   Suppose now that $\theta_1^*:=\Phi(f)$ is non-trivial and that $\nu$ is the $\theta$-Nielsen class rel $C$ that contains $z$.  
 Proposition~\ref{nonzero nielsen class} implies that the Nielsen class  of $f$ rel $C$ that corresponds to $\nu$  is non-trivial and compactly covered and hence that $\Fix(\ti f)$ is non-trivial, compact and $\ti \F$-invariant.  Theorem~\ref{thm: partial plane} implies that $\Fix(\ti \F) \ne \emptyset$.  Choose $\ti x \in \Fix(\ti \F)$ and let $x \in \Fix(\F)$ be its image in $M$.
\endproof

\begin{lemma}\label{invariant disk}  Suppose that  $\F$ is a finitely generated abelian
subgroup  of $\Homeo_+(S^2)$, that $X\subset S^2$ is a compact $\F$-invariant set and that $\Gamma \subset M:=S^2\setminus X$ is an essential simple closed curve.   If every element of $\F$ is isotopic rel $X$ to a homeomorphism $\theta$ that preserves $\Gamma$ and does not interchange the two complementary components of $\Gamma$  then $\Fix(\F) \ne \emptyset$.
\end{lemma}

\proof Let $D$ be one of the disks bounded by $\Gamma$ and let $M'$ be
the component containing $\Gamma$ of the surface obtained from $M$ by filling in the punctures
corresponding to the compact set $X_D :=X \cap D$.  Thus $D \subset
M'$ and $X_D$ is $\F$-invariant. Let $\ti M'$ be the universal
covering space of $M'$, let $\ti D \subset \ti M'$ be a lift of $D$
and let $\ti X_D$ be the lift of $X_D$ that is contained in $\ti D$.
Note that $\ti X_D$ is homeomorphic to $X_D$ and is therefore compact.
For each $f \in \F$ there is a unique lift $\ti f : \ti M' \to \ti M'$
that setwise fixes $\ti X_D$.  If $\F = \langle f_1, \ldots,
f_{n}\rangle$ then the uniqueness property implies that $\ti \F =
\langle \ti f_1, \ldots,\ti f_{n}\rangle$ is abelian.
Theorem~\ref{thm: partial plane} implies that $\Fix(\ti \F)$ and hence
$\Fix(\F)$ is non-empty.  \endproof

The following lemma is an application of the well known fact that any action of a finitely generated abelian group on a tree has a fixed point.

\begin{lemma} \label{tree argument} Suppose that $M$ is a genus zero surface, that $\F$ is a finitely generated abelian
subgroup  of $\Homeo_+(M)$ and that $R$ is a finite collection of simple closed disjoint curves in $M$ that is $\F$-invariant up to isotopy.  Let $\A$ be the union of disjoint open annulus neighborhoods of the elements of $R$.  Then there exists   either an element  $\Gamma$ of $R$ or a component $M_i$ of $M \setminus \A$ that is $\F$-invariant up to isotopy;   in the latter case, we may choose $M_i$ so that there is at most one component of $\partial M_i$ that is  $\F$-invariant up to isotopy.  
\end{lemma}

\proof   Let $T$ be the  finite tree that is dual to  $R$.  More precisely, $T$ has one vertex for each component $M_i$ of $M \setminus \A$ and an edge connecting the vertices corresponding to $M_i$ and $M_j$ if there exists $\Gamma \in R$ whose annulus neighborhood  has one boundary component in $M_i$ and the other in $M_j$.  If each edge of $T$ is assigned length one, then  $\F$ induces an action on $T$ by isometries.   

   Choose generators $\{f_i\}$ of $\F$. The fixed point set for the action of $f_1$ on $T$ is a non-empty subtree $T_1$. Since $\F$ is abelian, $T_1$ is $f_2$-invariant and   the fixed point set for the restricted action of $\F$ on $T_1$ is a non-empty subtree $T_2$ of $T_1$.  Repeating this argument shows that the set of points in $T$ that are fixed by  $\F$ is a non-empty tree $T^*$.

If $T^*$ is a single point and is not a vertex then it is the midpoint
of some edge of $T$.  In this case there is an $\F$-invariant $\Gamma
\in R$.  Otherwise, $T^*$ is either a single vertex or is a union of
edges of $T$.  Let $v$ be a vertex of $T^*$ whose valence in $T^*$ is
zero or one.  The subsurface $M_i$ corresponding to $v$ is
$\F$-invariant up to isotopy and components of $\partial M_i$
that are $\F$-invariant up to isotopy are in one-to-one correspondence
with edges of $T^*$ that are incident to $v$.  
\endproof

The index two part of Theorem~\ref{thm: main sphere} is handled by our next lemma.

\begin{lemma}  \label{index two}  If  a finitely generated abelian subgroup $\F$  of $\Homeo_+(S^2)$  has a finite index subgroup $F'$ such that $\Fix(\F') \ne \emptyset$ then it has such a subgroup with index at most two.  
\end{lemma}

\proof The $\F$-orbit $P$ of a point in $\Fix(\F')$ is finite.  If $P$
has cardinality one or two then its stabilizer has index at most two
and we are done.  We may therefore assume that $M := S^2 \setminus P$
has negative Euler characteristic.  Assume the notation of
Proposition~\ref{normal form} with $K = L = \emptyset$. (The existence
of normal forms in this case follows directly from the Thurston
classification theorem without reference to Proposition~\ref{normal
form} but it is convenient to use our usual notation.)  Denote
$\langle \theta_1,\dots,\theta_n\rangle$ by $\Theta$.

   If there is a $\Theta$-invariant element $\Gamma \in R$ then
Lemma~\ref{invariant disk} allows us to choose $\F'$ with index at
most two.  We are therefore reduced, by Lemma~\ref{tree argument}, to
the case that there is a $\Theta$-invariant subsurface $M_i$ whose
boundary does not have any $\Theta$-invariant components.  By
Lemma~\ref{moving on} there is a subgroup $\Theta'$ of $\Theta$ with
index at most two that fixes either a puncture, a component of
$\partial M$ or a point in the interior of $M_i$.  It therefore
suffices to show that the corresponding subgroup $\F'$ of $\F$ has a
global fixed point.  This is obvious if $\Theta'$ fixes a puncture and
follows from Lemma~\ref{invariant disk} if $\Theta'$ fixes a component
of $\partial M_i$.  In the remaining case, $\Theta'$ fixes a point $x$
in the interior of $M_i$ and Proposition~\ref{interior fixed point}
completes the proof.  \endproof

\noindent{\bf Proof of Theorem~\ref{thm: main sphere}}   We assume at first that $\F$ is finitely generated and argue by induction on the rank $n$ of $\F$.  The $n=1$ case is obvious so we assume inductively that the statement holds for rank $n-1$.  
 Choose generators $f_1,\ldots,f_n$ for $\F$ and denote $\langle f_1,\ldots ,f_{n-1}\rangle$ by $\F_{n-1}$.  By the inductive hypothesis there is a  subgroup $\F'_{n-1}$ of $\F_{n-1}$ with index at most two such that $K = \Fix(\F'_{n-1})$ is non-empty.  If $w(f,g) = 0$ for all $f, g \in \F$ then we may assume that $\F'_{n-1} = \F_{n-1}$.  Let $L = \Fix(f_n)$ and let $\F^*$ be the subgroup of $\F$ generated by $\F_{n-1}'$ and $f_n$.  Thus $\F^*$ has index at most two and equals $\F$   if $w(f,g) = 0$ for all $f, g \in \F$.    

If $K \cap L \ne \emptyset$ then $\F' = \F^*$ satisfies the conclusions of the theorem and we are done.  It therefore suffices to assume that $K \cap L = \emptyset$ and argue to a contradiction.   Assume the notation of Proposition~\ref{normal form} applied to $K,L$ and $\F^*$.  Let $\Theta'_{n-1}$ be the subgroup of normal forms corresponding to $\F'_{n-1}$, let $\theta_n$ correspond to $f_n$ and let $\Theta^*$ be the subgroup generated by $\Theta_{n-1}'$ and $\theta_n$.  

 If $R = \emptyset$ then $K \cup L$ is finite and its stabilizer has finite index in $\F^*$.   If $R \ne \emptyset$ then there is a finite index subgroup of $\Theta^*$ that setwise fixes each $M_i$.  Lemma~\ref{invariant disk} implies that the corresponding subgroup of $\F^*$ has a global fixed point.  In either case, Lemma~\ref{index two} implies that some subgroup of $\F$ with index at most two has a global fixed point.  

To complete the proof in the case that $\F$ is finitely generated, we must show that $\Fix(\F) \ne \emptyset$ if
and only if $w(f,g) = 0$ for all $f, g \in \F$.  The ``only if'' part is a
consequence of Corollary~\ref{cor: fp=>w=0}.  Thus we assume that
$w(f,g) = 0$ for all $f, g \in \F$ and prove that $\F^* = \F$ has a
global fixed point.  Suppose at first that $\Gamma \in R$ is
$\Theta^*$-invariant.  Since $K$ is non-empty, elements of
$\Theta'_{n-1}$ do not interchange the two complementary components of
$\Gamma$. Similarly, $\theta_n$ does not interchange the two
complementary components of $\Gamma$ because $L$ is non-empty.  It
follows that no element of $\Theta^*$ interchanges the two
complementary components of $\Gamma$.  Lemma~\ref{invariant disk} then
completes the proof.

We are therefore reduced, by Lemma~\ref{tree argument},  to the case that  there is  a $\Theta^*$-invariant subsurface $M_i$  whose boundary does not have any $\Theta^*$-invariant components.  Since $\theta_n$ fixes either a puncture in $M_i$ or a component of $\partial M_i$, it cannot be that $\Theta^*_{n-1}$ is the identity.  Similarly $\theta_n$ can not be the identity. Thus $M_i$ has finite type.  
 Lemma~\ref{moving on} implies that $\Theta^*$ fixes either a puncture, a boundary component or a point in the interior of $M_i$.  The first contradicts the assumption that $K \cap L = \emptyset$, the second  contradicts the choice of $M_i$ and the third contradicts the fact that $L = \Fix(\theta_n|_{M_i})$.   This completes the proof when $\F$ is finitely generated.

We complete the proof by showing the result for finitely generated
$\F$ implies the general result.  This follows from 
Proposition~\ref{prop: inf_gen} below applied with $n=1$ and $n=2$.
\endproof

\begin{prop}\label{prop: inf_gen}
Suppose $\Lambda$ is a compact metric space and $\F$ is a subgroup of
$\Homeo(\Lambda)$ with the property that every finitely generated
subgroup $\F_0$ has a subgroup of index $\le n$ with a global fixed point.
Then $\F$ has a subgroup of index $\le n$ with a global fixed point.
\end{prop}
\begin{proof}
We give a proof by contradiction.  A point $x \in \Lambda$ is a fixed
point of a subgroup of $\F$ of index $\le n$ if and only if the $\F$
orbit of $x$ contains at most $n$ points.  Thus if $x$ is not a fixed
point for a subgroup of $\F$ of index $\le n$ there are elements $f_{0,x},
f_{1,x},\dots f_{n,x} \in \F$ such that $f_{0,x}(x), f_{1,x}(x),\dots
f_{n,x}(x)$ are all distinct.  This is an open condition, i.e. there
is a neighborhood $U_x$ of $x$ such that for all $y \in U_x$ the
points $f_{0,x}(y), f_{1,x}(y),\dots f_{n,x}(y)$ are all
distinct. Assume that the result we want to prove is false, i.e. there
is no point of $\Lambda$ which is fixed by an index $n$ subgroup of
$\F.$ Then the collection $\{U_x\}$ is an open cover of $\Lambda$
which has a finite subcover $\{U_j\}_{j=1}^m$.  We will denote by
$f_{0,j}, f_{1,j},\dots f_{n,j}$ the corresponding elements of $\F$.
Then for any $x \in \Lambda$ there is a $j$ such that $f_{0,j}(x),
f_{1,j}(x),\dots f_{n,j}(x)$ are all distinct.  But this is impossible
since it would imply that the group $\F_0$ generated by
$\{f_{i,j}\}$ has no point fixed by an subgroup of index $\le n$
which contradicts the hypothesis.
\end{proof}

\noindent{\bf Proof of Corollary~\ref{cor: disk}} By
Proposition~\ref{prop: inf_gen} it suffices to prove the result when
$\F$ is finitely generated.  Let $\G \subset \Homeo_+(S^1)$ be the
group of restrictions to $S^1 = \partial \D^2$ of elements of $\F$.  If
every element of $\G$ has a fixed point then $\Fix(\G) \ne \emptyset.$
This is because the group $\G$ is abelian and hence there is a
$\G-$invariant measure $\mu$ on $S^1$.  But for each
$g \in \G \  \Fix(g) \ne \emptyset$ which implies $supp(\mu) \subset
\Fix(g).$  Hence $supp(\mu) \subset \Fix(\G).$

We are left with the case that there exists $f_0 \in \F$ such that
$\Fix(f_0) \cap \partial \D^2 = \emptyset$.  In that case we apply
Theorem~\ref{thm: partial plane} to $\D^2 \setminus \partial \D^2$ 
observing that $\Fix(f_0)$ is a compact $\F$ invariant set.
\endproof

\end{document}